\documentclass[11pt]{article}
\usepackage[a4paper, total={6in, 8in}]{geometry}
\usepackage[affil-it]{authblk}
\usepackage{graphicx} 
\usepackage{amsmath, amssymb, amsfonts, amsthm}
\usepackage[sorting=none]{biblatex}
\usepackage{enumerate}
\usepackage{xurl}
\usepackage[hyperindex, breaklinks]{hyperref}
\usepackage{tikz}
\usepackage{comment}
\usepackage{tikz-cd}
\usepackage{float}
\usepackage{nicematrix}
\hypersetup{
breaklinks=true,
colorlinks=true,
linkcolor=blue,
filecolor=magenta,
urlcolor=cyan,
pdfpagemode=FullScreen
}

\theoremstyle{plain}
\newtheorem{thm}{Theorem}[section]
\newtheorem*{thm*}{Theorem}
\newtheorem{lem}[thm]{Lemma}
\newtheorem{prop}[thm]{Proposition}
\newtheorem*{prop*}{Proposition}
\newtheorem{cor}[thm]{Corollary}

\newtheorem{conjecture}[thm]{Conjecture}

\theoremstyle{definition}
\newtheorem{defn}{Definition}[section]
\newtheorem{exmp}{Example}[section]
\newtheorem{notation}{Notation}[section]

\theoremstyle{remark}
\newtheorem*{rem}{Remark}

\title{On the $\Delta_a$ invariants in non-perturbative complex Chern-Simons theory}
\author{Shimal Harichurn \thanks{sharichurn.research@gmail.com}}

\date{}

\begin{document}

\maketitle

\begin{abstract}
Recently a set of $q$-series invariants, labelled by $\operatorname{Spin}^c$ structures, for weakly negative definite plumbed $3$-manifolds called the $\widehat{Z}_a$ invariants were discovered by Gukov, Pei, Putrov and Vafa. The leading rational power of the $\widehat{Z}_a$ invariants are invariants themselves denoted by $\Delta_a$. In this paper we further analyze the structure of these $\Delta_a$ invariants. We review some of the foundations of the $\Delta_a$ invariants and analyze their structure for a subclass of integer homology spheres. In particular, we provide a complete description of the $\Delta_0$ invariants for Brieskorn spheres. Along the way we show that the $\Delta_a$ invariants are not homology cobordism invariants, thereby answering an open question in the literature.
\end{abstract}

\tableofcontents

\section{Introduction}

The $\widehat{Z}_a$ invariants are a family of $q$-series invariants of weakly negative-definite plumbed $3$-manifolds associated to $\operatorname{Spin}^c$ structures that were introduced in \cite{bps-spectra} and take the form: 
\begin{equation}\label{base}
    \widehat{Z}_a(Y;q) = 2^{-\eta}q^{\Delta_a}(c_0 + c_1q^1 + c_2q^2 + \cdots) \in 2^{-\eta}q^{\Delta_a}\mathbb{Z}[[q]]
\end{equation} where $a$ is a $\operatorname{Spin}^c$ structure on a weakly negative-definite plumbed $3$-manifold $Y$, $\Delta_a \in \mathbb{Q}, \eta \in \mathbb{N} \cup \{0\}$ and $c_i \in \mathbb{Z}$. Both the numbers $\eta$ and $\Delta_a$ turn out to also be topological invariants. The $\widehat{Z}_a$ invariants can be seen as an extension of the Witten-Reshetikhin-Turaev (WRT) invariants \cite{witten_jones, rt_invariant}, $\tau_k$, away from roots of unity in the sense that if $b_1(Y) = 0$ for any $a \in T:= \operatorname{Spin}^c(Y)/\mathbb{Z}_2$, $\widehat{Z}_a(Y;q)$ converges in the unit disk $|q| < 1$ and for infinitely many $k$ the radial limits $\lim_{q \to r} \widehat{Z}_a(Y;q)$, where $r=e^{\frac{2\pi i}{k}}$, exist and the following (conjecturally) holds: $$\tau_k(Y) = \frac{1}{q^{1/2}-q^{-1/2}}\sum_{a \in T} \left(\sum_{b \in T}e^{2\pi i \ell k(a, a)}|\mathcal{W}_b|^{-1}S_{ab}\right)\widehat{Z}_a(Y;q)|_{q \to e^{2\pi i/k}}$$ where the coefficients $S_{ab}$ and $\mathcal{W}_x$ are given by $$S_{ab} = \frac{e^{2\pi i \ell k (a, b)} + e^{-2\pi i \ell k (a, b)}}{|\mathcal{W}_a| \cdot \sqrt{|H_1(Y;\mathbb{Z})|}} \ \ \ \text{ and } \ \ \ \mathcal{W}_x = \begin{cases}
2 \text{ if } x = \overline{x} \\
1 \text{ otherwise}
\end{cases}$$
where $\overline{x}$ denotes the $\operatorname{Spin}^c$ conjugate of $x$ (see \cite{murakami2024wittenreshetikhinturaevinvariantshomologicalblocks} for a proof of this for a large class of negative-definite plumbed manifolds). A further reason for studying the $\widehat{Z}_a$ invariants is that they are hoped to provide an extension of Khovanov homology \cite{khovanov1999categorificationjonespolynomial} to closed $3$-manifolds. Over the past few years the $\widehat{Z}_a$ invariants have been a very active area of research, see \cite{gukov2020twovariable, akhmechet2023latticecohomologyqseriesinvariants, Ekholm_2022, gukov2023hatzbplumbedmanifoldssplice, costin2023goingresurgentbridge, 3d_modularity, cheng20223manifoldsvoacharacters, Chun_2020, Ekholm_2022_Branches, quantum_z_g_invariants}.

In this paper, however, we will be mainly concerned with the $\Delta_a$ invariants. We will give a formal definition of these invariants later on, but, as \eqref{base} suggests, one can think of them as the unique rational pre-factored power of the $\widehat{Z}_a$ invariants. In \hyperref[read-delta]{Lemma \ref{read-delta}} we will formalize this notion further. 

\begin{exmp}\label{delta-lens-space}
As recorded in \cite[eq. 3.28]{cobordism}, if $Y = L(p, 1)$, one has $p$ many $\operatorname{Spin}^c$ structures of which $\widehat{Z}_a$ is only nonzero for two such $\operatorname{Spin}^c$ structures which we denote by $0$ and $1$. One has that $\widehat{Z}_0(Y;q) = -2q^{\frac{p-3}{4}}$ and $\widehat{Z}_1(Y;q) = q^{\frac{p-3}{4}}(2q^{\frac{1}{p}}).$ Thus $\Delta_0(Y) = \frac{p-3}{4}$ and $\Delta_1(Y) = \frac{p^2-3p +1}{4p}$.
\end{exmp}

\begin{exmp}\label{delta-s3}
By using the diffeomorphism between $L(1, 1)$ and $S^3$, and letting $0$ denote the trivial $\operatorname{Spin}^c$ structure on $S^3$ one obtains that $\widehat{Z}_0(S^3;q) = q^{-\frac{1}{2}}(2q-2)$ and hence $\Delta_0(S^3) = -\frac{1}{2}$.
\end{exmp}

As stated in \cite{cobordism}, the $\Delta_a$ invariants under the \textit{"3d Modularity Conjecture"} \cite{3d_modularity} are equal to the scaling dimension of associated log-VOA modules. Moreover, in \cite{cobordism} these $\Delta_a$ invariants were studied and related to the \textit{correction terms} $d$ in Heegaard Floer Homology. The correction terms are, in particular, homology cobordism invariants and the following relation was obtained (for an appropriate subclass of weakly negative-definite plumbed manifolds) \begin{equation}\label{main}
    \Delta_a(Y) = \frac{1}{2} - d(Y, a) \mod 1.
\end{equation} 
The relation in \eqref{main} is interesting in that it left open the possibility that $\Delta_a$ might itself be a homology cobordism invariant, which was a question asked in \cite{cobordism}. Searching for new homology cobordism invariants could potentially improve our understanding of the three-dimensional homology cobordism group $\Theta^3$ which is an area of active study. The structure of $\Theta^3$ has not yet been completely determined and any new progress in this direction would lead to a better understanding of low-dimensional topology. From what is known, the structure $\Theta^3$ is intimately connected to homology cobordism classes of Brieskorn spheres $\Sigma(b_1, b_2, b_3)$ defined by \[
\Sigma(b_1, b_2, b_3) := \{(x, y, z) \in \mathbb{C}^3 \mid x^{b_1} + y^{b_2} + z^{b_3} = 0\} \cap S^5.
\]
Indeed, as stated in \cite{_avk_2023} and originally proved in \cite{Furuta1990}, $\Theta^3$ has a $\mathbb{Z}^\infty$ summand generated by $\{\Sigma(2, 3, 6n-1)\}_{n=1}^\infty$ and it is an open question whether $\Theta^3$ is isomorphic to $\mathbb{Z}^\infty$. Thus it is natural when investigating a potential homology cobordism invariant to determine how to compute it on Brieskorn spheres. 

The goal of this paper is to review the structure of the $\Delta_a$ invariants and derive a formula for $\Delta_0$ on Brieskorn spheres (here $0$ denotes the trivial $\operatorname{Spin}^c$ structure) and use this to show that the $\Delta_a$ invariants are not homology cobordism invariants. The resulting formula for $\Delta_0$ on Brieskorn spheres will provide us with a playground of examples to study the relation \eqref{main} between the $\Delta_a$ invariants and the correction terms further. 

\subsection{Summary of results}

In what follows when $Y$ is an integer homology sphere  we will denote by $0$ the trivial $\operatorname{Spin}^c$ structure on $Y$ and write $\Delta_0$ and $\widehat{Z}_0$ for the corresponding invariants. In \hyperref[new-brieskorn-form]{Proposition \ref{new-brieskorn-form}}, we derive a formula for $\Delta_0$ for Brieskorn spheres which we partially present below.

\begin{prop*}
    Suppose $(b_1, b_2, b_3) \neq (2, 3, 5)$ where the integers $0 < b_1 < b_2 < b_3$ are pairwise relatively prime. Suppose the Brieskorn sphere $Y = \Sigma(b_1, b_2, b_3)$ has linking/framing matrix $M$ then we have that: 

$$\Delta_0= \frac{1}{4}\left(\sum_{i=1}^3 h_i  - 3s - \operatorname{Tr}(M) - \frac{b_2b_3}{b_1} - \frac{b_1b_3}{b_2} - \frac{b_1b_2}{b_3} \right) + \frac{\alpha_1^2}{4p}$$
wherein $p = b_1b_2b_3$, $\alpha_1 = b_1b_2b_3 - b_1b_2 - b_1b_3 - b_2b_3$, $h_i$ is the absolute value of the determinant of the linking matrix of the graph obtained by deleting a terminal vertex on the $i$-th leg from the decomposition $\Gamma$ of $\Sigma(b_1, b_2, b_3)$, and $s$ is the number of vertices in $\Gamma$. 
\end{prop*}

In the process we obtain a formula for $\widehat{Z}_0$ for Brieskorn spheres which is immediately of the form $q^{\Delta_0}\mathbb{Z}[[q]]$. In Examples \ref{sigma-2-9-11-inv} and \ref{sigma-3-7-8-inv} we compute $\widehat{Z}_0$ and $\Delta_0$ for $\Sigma(2, 9, 11)$ and $\Sigma(3, 7, 8)$.

In \hyperref[z-hat-cob-not]{Theorem \ref{z-hat-cob-not}} we show that neither $\widehat{Z}_a$ nor $\Delta_a$ are homology cobordism or cobordism invariants by providing some counterexamples. 

\begin{thm*}
$\widehat{Z}_a$ and $\Delta_a$ are neither homology cobordism invariants, nor cobordism invariants.
\end{thm*}

This answers a question raised in \cite{cobordism}. We further conjecture, in \hyperref[spinc-hom-cob-not]{Conjecture \ref{spinc-hom-cob-not}}, that $\Delta_a$ is not a $\operatorname{Spin}^c$ homology cobordism invariant. 

Using properties of the correction terms $d$ we are then able to prove the following (simple) proposition.

\begin{prop*}
    Let $\Gamma$ be an almost rational graph and let $Y := Y(\Gamma)$ be a negative definite plumbed $3$-manifold which is also an integral homology sphere, then $$\Delta_0(Y) = \frac{1}{2} \mod 1.$$
\end{prop*}

Finally in \hyperref[correction-terms-comparison]{Appendix \ref{correction-terms-comparison}}, we compare $\Delta_0$ and the correction terms $d$ for a class of Brieskorn spheres and in \hyperref[program]{Appendix \ref{program}} we compute the $\widehat{Z}_0$ and $\Delta_0$ invariants for various Brieskorn spheres which have not appeared thus far in the literature. \\ 

\noindent \textbf{Organization of the paper.}
In Section 2 we begin with a quick recap of the $\widehat{Z}_a$ invariants before defining the $\Delta_a$ invariants and providing a characterization for them. 

In Section 3 we will analyze $\Delta_0$ and $\widehat{Z}_0$ further for Brieskorn spheres $Y:= \Sigma(b_1, b_2, b_3)$ of the form $(b_1, b_2, b_3) \neq (2, 3, 5)$. In particular we will derive a formula for finding $\Delta_0(Y)$ and also put forward a proposition which puts $\widehat{Z}_0(Y;q)$ immediately into the form $q^{\Delta_0}P(q)$ for some $P(q) \in \mathbb{Z}[[q]]$. To digest these results we also give a couple of examples, those being $\Sigma(2, 9, 11)$ and $\Sigma(3,7,8)$ for which we compute $\Delta_0$ and $\widehat{Z}_0$.

In Section 4 we will show that the $\Delta_a$ invariants are not homology cobordism invariants by providing a counterexample. 

In Section 5 we will discuss the relation between $\Delta_a$ and the correction terms to Heegaard Floer homology. In this section we will also show that $\Delta_0(Y) = \frac{1}{2} \mod 1$ when $Y=Y(\Gamma)$ is a negative definite plumbed manifold arising from an almost rational graph $\Gamma$ and also an integer homology sphere. We will also discuss a conjecture that $\Delta_a$ is not a $\operatorname{Spin}^c$ homology cobordism invariant. \\

\noindent\textbf{Notation \& Conventions.} Unless otherwise specified, within this paper $Y = Y(\Gamma)$  will always mean a weakly negative-definite plumbed manifold with $b_1(Y) = 0$ which arises from a weighted tree $\Gamma$ with $s$ vertices. We will denote by $M$, the $s \times s$ linking/framing matrix of $\Gamma$ defined by 
$$M_{ij} = \begin{cases}
    0 \ \  \text{ if vertices $v_i, v_j$ are not adjacent},\\
    1 \ \ \text{ if vertices $v_i, v_j$ are adjacent},\\
    w_i \  \text{ if $i=j$},
\end{cases}$$
where $w_i$ is the weight associated to the vertex $v_i$. Further we will write $\vec{\delta} = (\deg(v))_{v \in \operatorname{Vert}(\Gamma)}$ for the degree vector of $\Gamma$. We will use the notation $\vec{a}$ to denote a representative (in $2\mathbb{Z}^s + \vec{\delta}$) of the $\operatorname{Spin}^c$ structure $a$ on $Y$. Finally if $Y$ is an integer homology sphere then we will denote by $0$ the trivial $\operatorname{Spin}^c$ structure on $Y$ and write $\Delta_0$ and $\widehat{Z}_0$ for the corresponding invariants. \\

\noindent \textbf{Acknowledgements.} The author would like to thank Sergei Gukov for guidance throughout the course of this project and feedback on a draft of this paper. The author would further like to thank Eveliina Peltola for helpful comments during the author's Master's thesis, from which this paper is derived. The author would like to thank Mrunmay Jagadale for stimulating discussions. Finally, the author would also like to thank Josef Svoboda for many technical discussions regarding the $\Delta_a$ invariants, in particular for making the author aware for the need of a reduced indexing set in the definition of the $\Delta_a$ invariants as well as for improvements to the direction of the paper. The author was supported by the FirstRand FNB 2020 Fund Education scholarship. \\

\section{The \texorpdfstring{$\Delta_a$}{Delta} invariant}

 We will begin with a quick recap of the $\widehat{Z}_a$-invariants. First let us recall the definitions of negative definite and weakly negative definite plumbed manifolds. The reader may find the definition of a plumbed manifold in \cite[p. 18]{gukov2020twovariable}.

 \begin{defn}
A weighted tree $\Gamma$, and its resulting plumbed manifold $Y(\Gamma)$ is called
\begin{enumerate}[(i)]
    \item \textit{negative definite} if the linking matrix $M$ it produces is negative definite.
    \item \textit{weakly negative definite} if the linking matrix $M$ it produces is invertible and $M^{-1}$ is negative definite on the subgroup of $\mathbb{Z}^s$ generated by the vertices of degree $\geq 3$. By this we mean that if $v_1, \dots, v_s$ are all the vertices of $\Gamma$ and $v_1, \dots, v_k$ where $k < s$ are all the vertices of degree $\geq 3$ then we require $M^{-1}$ is negative definite on the subgroup $\widetilde{\mathbb{Z}^k} := \{(l_{v_1}, \dots, l_{v_k}, 0, \dots, 0) \in \mathbb{Z}^s \mid l_{v_i} \in \mathbb{Z} \text{ for } 1 \leq i \leq k\}$ that is for every $\Vec{\ell} \in \widetilde{\mathbb{Z}^k}$ we require that $(\Vec{\ell}, M^{-1}\Vec{\ell}) < 0$.
\end{enumerate} 
\end{defn}

 The following is an adapted definition from \cite[p. 21]{gukov2020twovariable}.

\begin{defn}[Principal value]
Let $P(z)$ be a Laurent series which is holomorphic in some open set containing $\{z \in \mathbb{C} \mid 1 -c \leq |z| \leq 1 + c\}$ except possibly at $|z| =1$ (that is $P(z)$ can possibly be singular for $|z|=1$). We define $$\operatorname{v.p.} \int_{|z|=1} P(z) dz := \frac{1}{2} \left[\int_{|z|=1-\epsilon} P(z) dz + \int_{|z|=1+\epsilon} P(z) dz\right]$$ where $c>\epsilon > 0$ is some real number.
\end{defn}

\begin{rem}
The reason why the above definition doesn't depend on $\epsilon$ is by virtue of Cauchy's Theorem. Furthermore, if $P(z)$ is not singular for $|z| = 1$ then one sees that $
    \operatorname{v.p.} \int_{|z|=1} P(z) dz = \int_{|z|=1} P(z) dz $
\end{rem}

\begin{notation}
   If $f(z_0, \dots, z_s)$ is a complex function, then we will write $$\operatorname{v.p.} \int_{|z_0|=1} \cdots \int_{|z_s|=1}  f(z_0, \dots, z_s) dz_s \cdots dz_0$$ to mean $$\operatorname{v.p.} \int_{|z_0|=1} \operatorname{v.p.} \int_{|z_1|=1} \cdots \operatorname{v.p.} \int_{|z_s|=1}  f(z_0, \dots, z_s) dz_s \cdots dz_0.$$ 
\end{notation}

Now let us turn to the definition of the $\widehat{Z}_a$-invariants. These invariants are constructed in \cite[Appendix A]{bps-spectra}. The following definition has been taken with minor (but equivalent) modifications from \cite[p. 21]{gukov2020twovariable}.

\begin{defn}[The $\widehat{Z}_a$-invariants]\label{def-z-hat-invariants}
Let $\Gamma$ be a weighted tree with $s$ vertices which produces a weakly negative definite plumbed manifold $Y:= Y(\Gamma)$ with associated linking matrix $M$ such that $b_1(Y) = 0$.  Recall that   $\operatorname{Spin}^c(Y) \cong (2\mathbb{Z}^s + \Vec{\delta})/2M\mathbb{Z}^s$ (where $\Vec{\delta} = (\delta_v)_{v \in \operatorname{Vert}(\Gamma)} \in \mathbb{Z}^s$ denotes the vector comprised of the degrees of the vertices from $\Gamma$). For any $a \in \operatorname{Spin}^c(Y)$ let $\Vec{a} \in 2\mathbb{Z}^s + \Vec{\delta}$ be a representative of $a$, then letting $q$ be a formal variable, we define $\widehat{Z}_a(Y;q)$ by \begin{align*}\label{z-inv-def}
    \widehat{Z}_a(Y;q) &= (-1)^\pi q^{\frac{3\sigma - \operatorname{Tr}(M)}{4}}\\
    &\times \operatorname{v.p.}\oint_{|z_{v_1}|=1} \cdots \oint_{|z_{v_s}|=1} \left[\prod_{i=1}^s   \left(z_{v_i} - \frac{1}{z_{v_i}}\right)^{2-\deg(v_i)}\right] \cdot \Theta_{a}^{-M}(\Vec{z}) \frac{dz_{v_s}}{2\pi i z_{v_s}} \cdots \frac{dz_{v_1}}{2\pi i z_{v_1}} 
\end{align*}where $$\Theta_{a}^{-M}(\Vec{z}) = \sum_{\substack{\Vec{l} = (l_{v_1}, \dots, l_{v_s}), \\ \Vec{l} \in 2M\mathbb{Z}^s + \Vec{a}}}q^{-\frac{(\Vec{l}, M^{-1}\Vec{l})}{4}}\prod_{i=1}^sz_{v_i}^{l_{v_i}},$$ where in the above we have that  $\operatorname{v.p.}$ denotes taking the principal part of the integral, $\sigma$ is the signature of the matrix $M$ (the number of positive eigenvalues of $M$ minus the number of negative eigenvalues of $M$) and $\pi$ is the number of positive eigenvalues of $M$. 
\end{defn}

Let us now state a lemma which summarizes a few structural results on the $\widehat{Z}_a$ invariants. The second of which formalizes a comment made in \cite[p. 22]{gukov2020twovariable}.

\begin{lem}\label{z-hat-expanded-lemma}\label{finite-values-l}\label{wneg-def-delta}
Let $Y:= Y(\Gamma)$ be a weakly negative definite plumbed $3$-manifold with $b_1(Y) = 0$. Let $\vec{a}$ be a representative of $a \in \operatorname{Spin}^c(Y)$, then

\begin{enumerate}[(i)]
    \item $\widehat{Z}_a(Y;q)$ can be written in the form  \begin{equation}\label{z-hat-expanded}
    \widehat{Z}_a(Y;q) =  (-1)^\pi 2^{-\eta}q^{\frac{3\sigma - \operatorname{Tr}(M)}{4}}\sum_{\Vec{\ell} \in 2M\mathbb{Z}^s + \vec{a}} c_{\Vec{\ell}}\cdot q^{-\frac{(\Vec{\ell}, M^{-1}\Vec{\ell})}{4}},
\end{equation}
where $\eta \in \mathbb{N} \cup \{0\}$ and $c_{\vec{\ell}} \in \mathbb{Z}$.
\item  If $v_i$ is a vertex of $\Gamma$ which has $\deg(v_i) \leq 2$, then the number of values of $l_{v_i}$ in the tuple $\vec{\ell} = (l_{v_1}, \dots, l_{v_s}) \in 2M\mathbb{Z}^s + \vec{a}$ that yields a non-zero contribution to the $q$-series $\widehat{Z}_a(Y;q)$ is finite.
\end{enumerate}
    
\end{lem}

\begin{proof}
    (i) One applies the result (cf. \cite[p. 98]{lang_2010}) that $\int_{|z|=1}z^kdz = 0$ if $k \neq 1$ and it equals $2\pi i$ if $k = -1$  to the definition of the $\widehat{Z}_a$ invariants and the result follows. \\ \\
    (ii)  The proof we outlines builds upon a basic argument laid out in \cite[p. 22]{gukov2020twovariable}. Proving (ii) basically boils down to the fact that if we let $q$ be a formal variable, then an integral of the form
\begin{equation}\label{baby-integral}
    \oint_{|z|=1}\left(z-\frac{1}{z}\right)^m \sum_{l \in \mathbb{Z}} z^lq^l \frac{dz}{2\pi i z}
\end{equation} for $0 \leq m \leq 2$ is non-zero for only finitely many values of $l$. One then applies this general fact to the integrals which occur in $\widehat{Z}_a(Y;q)$ to complete the proof of the lemma.
\end{proof}

We now introduce a definition which will be useful in defining the $\Delta_a$ invariants to follow.

\begin{defn}[Reduced Indexing Set/Form] Suppose $\Gamma$ is a weighted tree with $s$ vertices and $Y:=Y(\Gamma)$ is a weakly negative definite plumbed manifold with $b_1(Y)=0$. Let $\vec{a}$ be a representative of a $\operatorname{Spin}^c$ structure $a$, we define a set $I \subseteq 2M\mathbb{Z}^s + \vec{a}$ to be a \textit{reduced indexing set} for $\widehat{Z}_a(Y;q)$ if we can write
    \begin{equation}\label{z_hat_reduced}
    \widehat{Z}_a(Y;q) =  \sum_{\vec{\ell} \in I} c_{\vec{\ell}}\cdot q^{\frac{3\sigma- \operatorname{Tr}(M)}{4} - \frac{(\vec{\ell}, M^{-1}\vec{\ell})}{4}}
    \end{equation}
    where $c_{\vec{\ell}}  \in \mathbb{Q} \setminus \{0\}$ for each $\vec{\ell} \in I$ and each $q$ power appearing in \eqref{z_hat_reduced} is unique. We say that $\widehat{Z}_a(Y;q)$ is written in \textit{reduced form} in this case.
\end{defn}

\begin{rem}
    One can see from Lemma \ref{z-hat-expanded-lemma} part (i) that a reduced indexing set $I$ for $\widehat{Z}_a(Y;q)$ always exists and that  $\widehat{Z}_a(Y;q)$ can always be put into reduced form. The reader should note, however, that a reduced indexing set  may not be unique.
\end{rem}

\begin{lem}\label{max_I_lemma}
    Let $Y:= Y(\Gamma)$ be a weakly negative definite plumbed $3$-manifold with $b_1(Y) = 0$. Let $\vec{a}$ be a representative of $a \in \operatorname{Spin}^c(Y)$ and suppose that $\widehat{Z}_a(Y;q)$ is given in reduced form as per \eqref{z_hat_reduced}. Then $\max_{\Vec{l} \in I} (\Vec{l}, M^{-1}\Vec{l})$ exists.
\end{lem}

\begin{proof}
    The proof we provide fleshes out the details of the basic argument laid out in \cite[p. 22]{gukov2020twovariable}. Suppose without loss of generality that $\Gamma$ has $k < s$ vertices of degree $\geq 3$ and if $v_1, \dots, v_s$ are the vertices of $\Gamma$ then $v_1, \dots, v_k$ are the vertices of degree $\geq 3$. Since $Y$ is a weakly negative definite plumbed manifold we have that $M^{-1}$ is negative definite on the subgroup $\widetilde{\mathbb{Z}^k} := \{(l_{v_1}, \dots, l_{v_k}, 0, \dots, 0) \in \mathbb{Z}^s \mid l_{v_i} \in \mathbb{Z} \text{ for } 1 \leq i \leq k\}$ that is for every $\Vec{\ell} \in \widetilde{\mathbb{Z}^k}$ we have that $(\Vec{\ell}, M^{-1}\Vec{\ell}) < 0$. In particular, negative-definiteness, shows that $\max_{\Vec{l} \in  \widetilde{\mathbb{Z}^k}} (\Vec{l}, M^{-1}\Vec{l})$ exists and hence that $\max_{\Vec{l} \in I \cap \widetilde{\mathbb{Z}^k}} (\Vec{l}, M^{-1}\Vec{l})$ exists. Now for any $\Vec{l} = (l_{v_1}, \dots,l_{v_k}, l_{v_{k+1}}, \dots l_{v_s}) \in I \subseteq 2M\mathbb{Z}^s + \Vec{a}$ we see from Lemma \ref{z-hat-expanded-lemma}, part (ii) that each of the $l_{v_{k+1}}, \dots, l_{v_s}$ can only take on finitely many different values from $\mathbb{Z}$. Thus one can view the indexing set $I$ as $I = \widetilde{\mathbb{Z}^k} \cup F$ where $F$ is some finite set. Thus we can conclude that $\max_{\Vec{l} \in I} (\Vec{l}, M^{-1}\Vec{l})$ exists from the fact that $\max_{\Vec{l} \in I \cap \widetilde{\mathbb{Z}^k}} (\Vec{l}, M^{-1}\Vec{l})$ exists.
\end{proof}




\subsection{The definition of the \texorpdfstring{$\Delta_a$}{Delta} invariant}

\begin{defn}[The $\Delta_a$-invariant]
Let $\Gamma$ be a weighted tree, with $s$ vertices which produces a weakly negative definite plumbed manifold $Y:=Y(\Gamma)$ with linking matrix $M$ and such that $b_1(Y) = 0$. Let $a \in \operatorname{Spin}^c(Y)$, and $\vec{a} \in 2\mathbb{Z}^s + \Vec{\delta}$ be a representative of $a$, suppose that $I$ is a reduced indexing set for $\widehat{Z}_a(Y;q)$, so that $\widehat{Z}_a(Y;q)$ is given in reduced form as per \eqref{z_hat_reduced}, then we define $$\Delta_a = \frac{3\sigma - \operatorname{Tr}(M)}{4} - \max_{\Vec{l} \in I} \frac{(\Vec{l}, M^{-1}\Vec{l})}{4}$$ where $(\Vec{l}, M^{-1}\Vec{l})$ means the inner product of $\Vec{l}$ with $M^{-1}\Vec{l}$.
\end{defn}

\begin{rem}
\begin{enumerate}[a.)]
    \item We will sometimes write $\Delta_a(Y)$ to show the explicit dependence on the manifold $Y$, but we will for the most part simply write $\Delta_a$ instead.
    \item The definition we provided above was produced in a slightly different manner in \cite[Appendix A]{bps-spectra}. In \cite[Appendix A]{bps-spectra} the indexing set $I$ was not mentioned in the maximum $\max_{\Vec{l} \in I} \frac{(\Vec{l}, M^{-1}\Vec{l})}{4}$ which occurs in the definition of $\Delta_a$, though we believe that this indexing set was implicitly assumed, as is done in \cite{cobordism}.
    \item Lemma \ref{max_I_lemma} ensures that the maximum used in the definition above exists.
    \item Note that the maximum used in the definition of $\Delta_a$ is independent of the choice of reduced indexing set and so the definition is well-defined.
\end{enumerate}

\end{rem}

The following lemma shows that the $\Delta_a$-invariants really are invariants of weakly negative definite plumbed manifolds.

\begin{lem}
Let $\Gamma_0$, $\Gamma_1$ be two weighted trees, which produce weakly negative definite plumbed manifolds $Y_0 := Y(\Gamma_0)$ and $Y_1 := Y(\Gamma_1)$ with linking matrices $M_0$ and $M_1$ respectively. If $Y_0$ is diffeomorphic to $Y_1$ then we have that $\Delta_a(Y_0) = \Delta_a(Y_1)$
\end{lem}

\begin{proof}
This follows immediately since we have that $\widehat{Z}_a(Y_0;q) = \widehat{Z}_a(Y_1;q)$ because the fact that $Y_0$ is diffeomorphic to $Y_1$ implies that $\Gamma_0$ and $\Gamma_1$ are related by a sequence of Neumann moves and the $\widehat{Z}_a$-invariants are invariant under Neumann moves.
\end{proof}

It's well known that one can $\widehat{Z}_a(Y;q)$ into an algebraic object of the form $2^{-\eta}q^{\Delta_a}P(q)$ where $\Delta_a$ is some rational power, $\eta \in \mathbb{N} \cup \{0\}$ and $P(q) \in \mathbb{Z}[[q]]$. For the sake of completeness, given the subject matter of this paper, we provide a proof of this below. 

\begin{prop}\label{prop-z-hat-alg}
We have that $\widehat{Z}_a(Y;q) \in 2^{-\eta}q^{\Delta_a}\mathbb{Z}[[q]]$ where $\eta \in \mathbb{N} \cup \{0\}$ and $2^{-\eta}q^{\Delta_a}\mathbb{Z}[[q]] := \left\{2^{-\eta}q^{\Delta_a}P(q) \mid P(q) \in \mathbb{Z}[[q]]\right\}.$
\end{prop}

\begin{proof}
We can assume $\widehat{Z}_a(Y;q)$ is given in reduced form, that is
\begin{equation}\label{q-power}
       \widehat{Z}_a(Y;q) =  (-1)^\pi 2^{-\eta}q^{\frac{3\sigma - \operatorname{Tr}(M)}{4}}\sum_{\Vec{\ell} \in I} c_{\Vec{\ell}}\cdot q^{-\frac{(\Vec{\ell}, M^{-1}\Vec{\ell})}{4}} 
\end{equation}
where $I \subseteq 2M\mathbb{Z}^s + \Vec{a}$ is an indexing set such that $c_{\vec{\ell}} \neq 0$  for each $\vec{\ell} \in I$, $\Vec{a} \in 2\mathbb{Z}^s + \Vec{\delta}$ is a \textit{fixed} representative of $a$ and $\eta \in \mathbb{N} \cup \{0\}$. Since $Y$ is a weakly negative definite plumbed manifold by Lemma \ref{max_I_lemma} we see that $\gamma:= \max_{\Vec{\ell} \in I} (\Vec{\ell}, M^{-1}\Vec{\ell})$ exists. Thus we can factorize $\widehat{Z}_a(Y;q)$ into the form 
\begin{equation}\label{factorized}
   \widehat{Z}_a(Y;q) =  (-1)^\pi 2^{-\eta}q^{\frac{3\sigma - \operatorname{Tr}(M) - \gamma}{4}}\sum_{\Vec{\ell} \in I} c_{\Vec{\ell}}\cdot q^{\frac{\gamma - (\Vec{\ell}, M^{-1}\Vec{\ell})}{4}} 
\end{equation}
Now let  $\vec{\ell}_0 \in I$ be such that $\gamma = (\Vec{\ell_0}, M^{-1}\Vec{\ell_0})$. Then to show that $\widehat{Z}_a(Y;q) \in 2^{-\eta}q^{\Delta_a}\mathbb{Z}[[q]],$ where $\Delta_a = \frac{3\sigma - \operatorname{Tr}(M) - \gamma}{4}$ it suffices to show that for any $\vec{\ell} \in I$ that $\gamma-(\vec{\ell}, M^{-1}\vec{\ell}) = (\Vec{\ell_0}, M^{-1}\Vec{\ell_0}) - (\vec{\ell}, M^{-1}\vec{\ell})   \in 4 \mathbb{Z}_+$ where $\mathbb{Z}_+ = \{x \in \mathbb{Z} \mid x \geq 0\}$ because then it would follow that all the exponents of $q$ on the right hand side of equation \eqref{factorized} are non-negative integers.\\ \\ To that end let $\vec{\ell} \in I$ be given arbitrarily. First notice that since $\gamma:= \max_{\Vec{\ell} \in I} (\Vec{\ell}, M^{-1}\Vec{\ell})$ we have $\gamma-(\vec{\ell}, M^{-1}\vec{\ell}) \geq 0$. Then secondly notice that we can express $\vec{\ell} =  \Vec{\ell_{0}} + 2Mv$ where $v \in \mathbb{Z}^s$. The reason for this is that both $\vec{\ell}$ and $\Vec{\ell_{0}}$ are elements of the set $2M\mathbb{Z}^s + \Vec{a}$ and since $\Vec{a}$ is fixed we have that $\vec{\ell} - \Vec{\ell_{0}} \in 2M\mathbb{Z}^s$, which implies that $\vec{\ell} =  \Vec{\ell_{0}} + 2Mv$. One can then compute (using the fact that $M$ is a symmetric matrix and properties of the transpose) that  
\begin{align*}
    \left(\vec{\ell}, M^{-1}\vec{\ell}\right) &= \vec{\ell}^TM^{-1}\vec{\ell}
    = (\Vec{\ell_{0}} + 2Mv)^TM^{-1}(\Vec{\ell_{0}} + 2Mv) 
= \gamma + 4\Vec{\ell_{0}}^Tv  + 4(v, Mv).
\end{align*}
Noting that $\Vec{\ell_{0}}^Tv \in \mathbb{Z}$ and $(v, Mv) \in \mathbb{Z}$, we can thus see that $\gamma-(\vec{\ell}, M^{-1}\vec{\ell}) = -4\Vec{\ell_{0}}^Tv  - 4(v, Mv)  \in 4 \mathbb{Z}$. Combining this with the fact that $\gamma-(\vec{\ell}, M^{-1}\vec{\ell}) \geq 0$ implies that $\gamma-(\vec{\ell}, M^{-1}\vec{\ell}) \in 4\mathbb{Z}_+$. This completes the proof of the proposition.
\end{proof}

\subsection{A characterization of the \texorpdfstring{$\Delta_a$}{Delta} invariant}

The following lemma provides a very useful way to simply 'read off' the $\Delta_a$-invariant from a given expression of $\widehat{Z}_a(Y;q)$. A small statement within the statement of the lemma below appeared in \cite[Equation 1.2]{cobordism} in a different form.

\begin{lem}\label{read-delta}
Suppose that $\widehat{Z}_a(Y;q) = q^\delta(c_0 + c_1q^{n_1} + c_2q^{n_2} + \cdots )$ where $\delta \in \mathbb{Q}$, $c_i \in \mathbb{Q} \setminus \{0\}$ and $n_i \in \mathbb{N} \setminus \{0\}$ for each $i$, then $\delta = \Delta_a$.
\end{lem}

\begin{proof}
Without loss of generality, after simplifying the expression for $\widehat{Z}_a(Y;q)$ in the statement of the Lemma, we can further assume that $n_i \neq n_j$ for $i \neq j$. From $\widehat{Z}_a(Y;q) = q^\delta(c_0 + c_1q^{n_1} + c_2q^{n_2} + \cdots )$, simply multiply the $q^\delta$ factor out to get \begin{equation}\label{read-eq-1}
    \widehat{Z}_a(Y;q) = c_0q^{\delta} + c_1q^{n_1 + \delta} + c_2q^{n_2 + \delta} + \cdots = q^\delta \sum_{i \geq 1}c_iq^{n_i}.
\end{equation}
We see, from Lemma \ref{z-hat-expanded-lemma} part (i), that for each $i \geq 1$ we must have \begin{equation}\label{n-i}
    n_i + \delta = \frac{3\sigma -\operatorname{Tr}(M) - (\vec{\ell}_i, M^{-1}\vec{\ell}_i)}{4}
\end{equation} for some $\vec{\ell}_i \in 2M\mathbb{Z}^s + \vec{a}$. Further we see that \begin{equation}\label{delta-equals}
    \delta = \frac{3\sigma -\operatorname{Tr}(M) - (\vec{\ell}_0, M^{-1}\vec{\ell}_0)}{4}
\end{equation} for some $\vec{\ell}_0 \in 2M\mathbb{Z}^s + \vec{a}$. Equation \eqref{delta-equals} then implies that $$n_i = n_i + \delta - \delta = \frac{(\vec{\ell}_0, M^{-1}\vec{\ell}_0) - (\vec{\ell}_i, M^{-1}\vec{\ell}_i)}{4}.$$ The assumption that $n_i \in \mathbb{N} \setminus \{0\}$ then implies that $n_i >0$ and hence that $(\vec{\ell}_0, M^{-1}\vec{\ell}_0)> (\vec{\ell}_i, M^{-1}\vec{\ell}_i)$ for each $i > 1$. Furthermore, since $n_i \neq n_j$ for $i \neq j$ and $c_i \neq 0$ for all $i$, the set $$I := \{\vec{\ell}_i\}_{i \geq 0}$$ is a reduced indexing set for $\widehat{Z}_a(Y;q)$. From the reasoning above we see that $(\vec{\ell}_0, M^{-1}\vec{\ell}_0) = \max_{\Vec{l} \in I } \frac{1}{4}(\Vec{l}, M^{-1}\Vec{l})$ and hence we have that $$\delta = \frac{3\sigma -\operatorname{Tr}(M)}{4} - \max_{\Vec{l} \in I } \frac{(\Vec{l}, M^{-1}\Vec{l})}{4} = \Delta_a$$ as desired. 
\end{proof}

\begin{rem}
With the proof of this lemma, our earlier remark in the introduction that one can view $\Delta_a$ as the unique rational pre-factored power of $q$ in $\widehat{Z}_a(Y;q)$ is justified.
\end{rem}

\section{\texorpdfstring{$\Delta_0$}{Delta} for Brieskorn Spheres}

Given relatively co-prime integers $0 < b_1 < b_2 < b_3$ one can define the \textit{Brieskorn sphere} $\Sigma(b_1, b_2, b_3) := \{(x, y, z) \in \mathbb{C}^3 \mid x^{b_1} + y^{b_2} + z^{b_3} = 0\} \cap S^5.$ Brieskorn spheres are negatively definite plumbed $3$-manifolds which are also integral homology spheres. Hence for Brieskorn spheres there is a unique $\operatorname{Spin}^c$ structure, which we label as $0$. Hence there is only a single corresponding $\widehat{Z}_a$-invariant that being $\widehat{Z}_0(\Sigma(b_1, b_2, b_3);q)$ and a single $\Delta_0$-invariant.  The method for finding a plumbing description for $\Sigma(b_1, b_2, b_3)$, was mentioned in \cite[p. 21]{gukov2020twovariable}. Let us explain this method. The Brieskorn sphere $\Sigma(b_1, b_2, b_3)$ is the Seifert manifold $M(b; \frac{a_1}{b_1},\frac{a_2}{b_2}, \frac{a_3}{b_3} )$ where $b < 0$ and $a_1, a_2, a_3 >0$ are \textit{chosen} to satisfy the equation
\begin{equation}
    \label{brieskorn-eq}
    b_1b_2b_3 \cdot b + b_2b_3a_1 + b_1b_3a_2 + b_1b_2a_3 = -1.
\end{equation}
The integers $b, a_1, a_2, a_3$ are chosen as there is far more than a single solution to equation \ref{brieskorn-eq} above (thus there is also more than one plumbing description). 

\begin{figure}[H]
\centering

\begin{tikzpicture}

        \node (0A) at (-4, -0.5) {$v_0$};
        \node (0B) at (-2, 1.5) {$v_{1_1}$};
        \node (0C) at (-2, -0.5) {$v_{2_1}$};
        \node (0D) at (-2, -2.5) {$v_{3_1}$};
        \node (0E) at (0, -0.5) {$v_{2_2}$};
        \node (0F) at (0, -2.5) {$v_{3_2}$};
        \node (0G) at (0, 1.5) {$v_{1_2}$};
        \node (0H) at (4, 1.5) {$v_{1_{s_1}}$};
        \node (0I) at (4, -2.5) {$v_{3_{s_3}}$};
        \node (0J) at (4, -0.5) {$v_{2_{s_2}}$};

        \node (0A) at (-4, 0.5) [text=teal]{$b$};
        \node (0B) at (-2, 2.5) [text=teal]{$-k_{1_1}$};
        \node (0C) at (-2, 0.5) [text=teal]{$-k_{2_1}$};
        \node (0D) at (-2, -1.5) [text=teal]{$-k_{3_1}$};
        \node (0E) at (0, 0.5) [text=teal]{$-k_{2_2}$};
        \node (0F) at (0, -1.5) [text=teal]{$-k_{3_2}$};
        \node (0G) at (0, 2.5) [text=teal]{$-k_{1_2}$};
        \node (0H) at (4, 2.5) [text=teal]{$-k_{1_{s_1}}$};
        \node (0I) at (4, -1.5) [text=teal]{$-k_{3_{s_3}}$};
        \node (0J) at (4, 0.5) [text=teal]{$-k_{2_{s_2}}$};
        
        \begin{scope}[every node/.style={circle, thick, draw, fill=black, inner sep =0pt, minimum size=5pt}]
        \node (1) at (-4, 0) {};
        \node (2A) at (-2, 2) {};
         \node (2B) at (-2, 0) {};
          \node (2C) at (-2, -2) {};
          \node (3A) at (0, 0) {};
          \node (3B) at (0, -2) {};
          \node (3C) at (0, 2) {};
        \end{scope}

         \begin{scope}[every node/.style={circle, thick, draw, fill=black, inner sep =0pt, minimum size=5pt}]
        \node (M1) at (2, 0) {};
        \node (M2) at (2, 2) {};
        \node (M3) at (2, -2) {};
        \end{scope}

          \begin{scope}[every node/.style={circle, thick, draw, fill=black, inner sep =0pt, minimum size=5pt}]
        \node (MM1) at (4, 0) {};
        \node (MM2) at (4, 2) {};
        \node (MM3) at (4, -2) {};
        \end{scope}

        \draw (1) -- (2A);
        \draw (1) -- (2B);
         \draw (1) -- (2C);
         \draw (2B) -- (3A);
           \draw (2C) -- (3B);
           \draw (2A) -- (3C);
           \draw[dashed] (3A) -- (MM1);
           \draw[dashed] (3B) -- (MM3);
           \draw[dashed] (3C) -- (MM2);
    \end{tikzpicture}

\caption{Plumbing description of the Brieskorn manifold $Y(\Gamma) := \Sigma(b_1, b_2, b_3) = M(b;\frac{a_1}{b_1},\frac{a_2}{b_2}, \frac{a_3}{b_3})$. We call the vertices which correspond to the weights $-k_{i_{s_i}}$ for $1 \leq i \leq 3$, the \textit{terminal vertices} of the graph $\Gamma$. Sometimes we also call it the 'terminal vertices of the plumbing description' of the Brieskorn sphere $\Sigma(b_1, b_2, b_3)$.}
\label{fig-brie-decomp}
\end{figure}
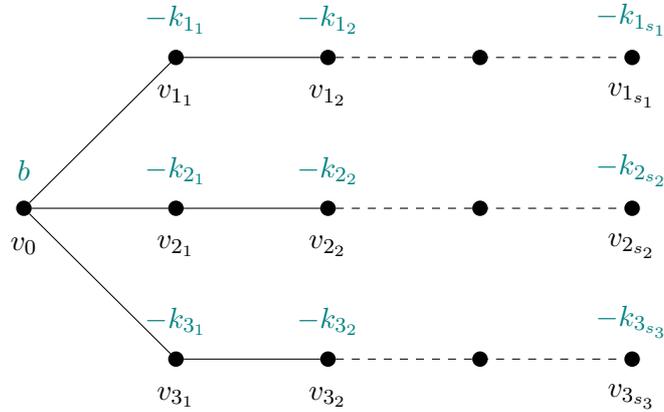

Then the plumbing graph $\Gamma$ (see Figure \ref{fig-brie-decomp}) which produces $\Sigma(b_2, b_2, b_3)$ has a central vertex labelled by $b$ and three "legs" whose vertices are given by the integers $-k_{i_1}, \dots, -k_{i_{s_{i}}}$ for $1 \leq i \leq 3$ which come from the continued fraction decomposition of $\frac{b_i}{a_i}$, i.e. $$\frac{b_i}{a_i} = [k_{i_1}, \dots, k_{i_{s_i}}] =k_{i_1} - \cfrac{1}{k_{i_2} - \cfrac{1}{\cdots - \cfrac{1}{k_{i_{s_i}}}}}.$$

Let us give an example to illustrate this procedure.

\begin{exmp}\label{sigma-2-9-11}
Consider the Brieskorn sphere $\Sigma(2, 9, 11)$. The integers $b = -1$, $a_1 = 1, a_2 = 2$ and $a_3 = 3$ satisfy equation \ref{brieskorn-eq} with $b_1 = 2, b_2 = 9$ and $b_3=11$. One can then compute the continued fractions $\frac{b_1}{a_1} = 2 = [2]$, $\frac{b_2}{a_2} = \frac{9}{2} = [5,2]$ and $\frac{b_3}{a_3} = \frac{11}{3} = [4,3]$ to produce the following plumbing graph for $\Sigma(2, 9, 11)$.

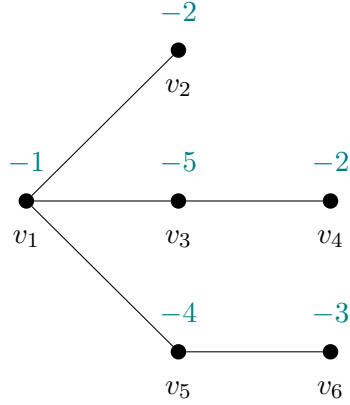
\begin{figure}[H]
\centering

    \begin{tikzpicture}

        \node (0A) at (-4, -0.5) {$v_1$};
        \node (0B) at (-2, 1.5) {$v_2$};
        \node (0C) at (-2, -0.5) {$v_3$};
        \node (0D) at (-2, -2.5) {$v_5$};
        \node (0E) at (0, -0.5) {$v_4$};
        \node (0F) at (0, -2.5) {$v_6$};

        \node (0A) at (-4, 0.5) [text=teal]{$-1$};
        \node (0B) at (-2, 2.5) [text=teal]{$-2$};
        \node (0C) at (-2, 0.5) [text=teal]{$-5$};
        \node (0D) at (-2, -1.5) [text=teal]{$-4$};
        \node (0E) at (0, 0.5) [text=teal]{$-2$};
        \node (0F) at (0, -1.5) [text=teal]{$-3$};
        
        \begin{scope}[every node/.style={circle, thick, draw, fill=black, inner sep =0pt, minimum size=5pt}]
        \node (1) at (-4, 0) {};
        \node (2A) at (-2, 2) {};
         \node (2B) at (-2, 0) {};
          \node (2C) at (-2, -2) {};
          \node (3A) at (0, 0) {};
          \node (3B) at (0, -2) {};
        \end{scope}

        \draw (1) -- (2A);
        \draw (1) -- (2B);
         \draw (1) -- (2C);
         \draw (2B) -- (3A);
           \draw (2C) -- (3B);
    \end{tikzpicture}

\caption{Plumbing graph for $\Sigma(2, 9, 11)$. The terminal vertices of $\Gamma$ are $v_2, v_4, v_6$.}
\end{figure}
\end{exmp}

\begin{notation}\label{brieskorn-notation}
For the rest of this section let us fix the following notation: suppose we are given co-prime integers $b_i$, $1 \leq i \leq 3$ which satisfy $0 < b_1 < b_2 < b_3$ and we are considering the Brieskorn sphere $\Sigma(b_1, b_2, b_3)$. We then define $p := b_1b_2b_3$ and \begin{align*}
    \alpha_1 &:= b_1b_2b_3 - b_1b_2 - b_1b_3 - b_2b_3, \\
    \alpha_2 &:= b_1b_2b_3 + b_1b_2 - b_1b_3 - b_2b_3, \\
    \alpha_3 &:= b_1b_2b_3 - b_1b_2 + b_1b_3 - b_2b_3, \\
    \alpha_4 &:= b_1b_2b_3 + b_1b_2 + b_1b_3 - b_2b_3. 
\end{align*}
\end{notation}

In \cite[Proposition 4.8]{gukov2020twovariable} a simplification of the general formula for $\widehat{Z}_0(\Sigma(b_1, b_2, b_3);q)$ was put forward which incorporated functions known as \textit{false theta functions}.\footnote{As mentioned in \cite{gukov2020twovariable}, another derivation of the formula for $\widehat{Z}_0(\Sigma(b_1, b_2, b_3);q)$ was given by Chung in \cite{chung2020bps}.} These functions are defined in \cite[Equations 53 and 54]{gukov2020twovariable} as
\begin{equation}\label{theta-func-1}
    \widetilde{\Psi}^{(a)}_p(q) := \sum_{n=0}^\infty \psi^{(a)}_{2p}(n) q^{\frac{n^2}{4p}} \in q^{\frac{a^2}{4p}}\mathbb{Z}[[q]]
\end{equation} where 
\begin{equation}\label{theta-func-2}
\psi^{(a)}_{2p}(n) = \begin{cases}
\ \ \  1   \text{ if } n= \ \ a +m\cdot 2p \text{ for } m \in \mathbb{Z}\\
\ -1 \text{ if } n= -a +m\cdot 2p \text{ for } m \in \mathbb{Z}\\
\ \ \ 0 \text{ otherwise}
\end{cases}    
\end{equation}
and, following the convention in \cite[Equation 55]{gukov2020twovariable}, the notation $\widetilde{\Psi}^{c_1(a_1)+ c_2(a_2) + \cdots}_p(q)$ is used as a shorthand for $c_1\widetilde{\Psi}^{(a_1)}_p(q) + c_2\widetilde{\Psi}^{(a_2)}_p(q) + \cdots$. We noticed however, that there was a sign error contained in the proof of  \cite[Proposition 4.8]{gukov2020twovariable}. Correcting this sign error we state the simplified formula below.

\begin{thm}\label{corrected-brie}
For $(b_1, b_2, b_3) \neq (2, 3, 5)$, consider the Brieskorn sphere $Y = \Sigma(b_1, b_2, b_3)$ which has linking matrix $M$, where the integers $0 < b_1 < b_2 < b_3$ are pairwise relatively prime. Then we have that \begin{equation}\label{z-hat-form-prop-4-8}
    \widehat{Z}_0(Y;q) = q^\xi \cdot \left(\widetilde{\Psi}_{b_1b_2b_3}^{(\alpha_1) - (\alpha_2) - (\alpha_3) + (\alpha_4)}(q)\right)
\end{equation} where $\alpha_i$ for $1 \leq i\leq 3$ is defined as in Notation \ref{brieskorn-notation} and
$$\xi = \frac{1}{4}\left(\sum_{i=1}^3 h_i  - 3s - \operatorname{Tr}(M) - \frac{b_2b_3}{b_1} - \frac{b_1b_3}{b_2} - \frac{b_1b_2}{b_3} \right) \in \mathbb{Q}$$ wherein $h_i$ is the absolute value of the determinant of the linking matrix of the graph obtained by deleting a terminal vertex on the $i$-th leg from the plumbing description $\Gamma$ which produces $\Sigma(b_1, b_2, b_3)$, and $s$ is the number of vertices in $\Gamma$. \footnote{Note that in \cite{gukov2020twovariable}, the notation $\Delta$ was adopted instead of the notation $\xi$ we use here. One should however not confuse the notation $\Delta$ used in \cite[Proposition 4.8]{gukov2020twovariable} with the $\Delta_0$ invariant for $\Sigma(b_1, b_2, b_3)$. From what we will show below one can see that $\Delta$ and $\Delta_0$ are related but not equal.}
\end{thm}

The sign error that occurred in \cite[Proposition 4.8]{gukov2020twovariable} occurs towards the end of the proof.

We do not consider $(b_1, b_2, b_3) = (2, 3, 5)$ because in that case the form of $\widehat{Z}_0(\Sigma(b_1, b_2, b_3);q)$ which appears in equation \eqref{z-hat-form-prop-4-8} is not entirely correct, namely there is an extra term that one needs to include. One can see this extra term in the statement of \cite[Proposition 4.8]{gukov2020twovariable}. To simplify our proofs below and for brevity we exclude this case for the rest of the section.

\subsection{A formula for \texorpdfstring{$\Delta_0$}{Delta} for Brieskorn spheres}

We first state a few useful results that will help us.

\begin{lem}\label{gamma-min}\label{alpha-i-geq-0}\label{brie-pow-int} 
    Suppose we are given co-prime integers $b_i$, $1 \leq i \leq 3$ which satisfy $0 < b_1 < b_2 < b_3$. Let $p$ and $\alpha_i$ for $1 \leq i \leq 4$ be given as in Notation \ref{brieskorn-notation}. Then
    \begin{enumerate}[(i)]
        \item $\min\{\alpha_1,\alpha_2,\alpha_3,\alpha_4\} = \alpha_1$.
        \item $\frac{\alpha_i^2 -\alpha_1^2}{4p} \in \mathbb{Z}$ for $1 \leq i \leq  4$. 
        \item Suppose further that $(b_1, b_2, b_3) \neq (2, 3, 5)$, then $\frac{1}{b_1} + \frac{1}{b_2} + \frac{1}{b_3} < 1$ and $0 < \alpha_i < 2p$ for $1 \leq i \leq 4$. 
    \end{enumerate}
\end{lem}

\begin{proof}
(i) We have the following series of equalities $$\alpha_1 = \alpha_2 - b_1b_2 = \alpha_3 - b_1b_3 = \alpha_4 -b_1b_2- b_1b_3.$$
Then the fact that $b_i > 0$ for $1 \leq i \leq 3$ implies that $\alpha_1 < \alpha_j$ for $2 \leq j \leq 4$ which completes the proof. \\ \\
(ii) Note that to show that $\frac{\alpha_i^2 -\alpha_1^2}{4p} \in \mathbb{Z}$ for $i = 1, 2, 3, 4$ it suffices to show that $4p | (\alpha_i^2 - \alpha_1^2)$ for $i = 2,3,4$. Since $\alpha_i^2 - \alpha_1^2 = (\alpha_i+\alpha_1)(\alpha_i-\alpha_1)$, if we can show that $4p$ divides $(\alpha_i+\alpha_1)(\alpha_i-\alpha_1)$ for $i = 2,3,4$ then we are done.  To that end observe that we have the following series of equalities:
\begin{align*}
    (\alpha_2+\alpha_1)(\alpha_2-\alpha_1) &= (2b_1b_2b_3-2b_1b_3-2b_2b_3)(2b_1b_2) \\
    (\alpha_3+\alpha_1)(\alpha_3-\alpha_1) &= (2b_1b_2b_3-2b_1b_2-2b_2b_3)(2b_1b_3)\\
    (\alpha_4+\alpha_1)(\alpha_4-\alpha_1) &= (2b_1b_2b_3-2b_2b_3)(2b_1b_2+2b_1b_3)
\end{align*}
In all the cases above one can just expand the right hand side and check directly that $4p$ where $p=b_1b_2b_3$ divides all the quantities above. \\ \\
(iii) Proof omitted - follows from elementary number theory arguments. 
\end{proof}

\begin{prop}\label{new-brieskorn-form}\label{closed-delta-brie}
Suppose $(b_1, b_2, b_3) \neq (2, 3, 5)$ where the integers $0 < b_1 < b_2 < b_3$ are pairwise relatively prime. Suppose the Brieskorn sphere $Y = \Sigma(b_1, b_2, b_3)$ has linking matrix $M$ then we have that: 

\begin{align*}
       \widehat{Z}_0(Y;q) &=  q^{\Delta_0}\cdot \left[q^{-\frac{\alpha_1^2}{4p}}\left(\widetilde{\Psi}_{b_1b_2b_3}^{(\alpha_1) - (\alpha_2) - (\alpha_3) + (\alpha_4)}(q)\right)\right]
\end{align*}
where $$\Delta_0= \xi + \frac{\alpha_1^2}{4p}$$
wherein $p = b_1b_2b_3$ and $\xi = \frac{1}{4}\left(\sum_{i=1}^3 h_i  - 3s - \operatorname{Tr}(M) - \frac{b_2b_3}{b_1} - \frac{b_1b_3}{b_2} - \frac{b_1b_2}{b_3} \right) \in \mathbb{Q}$
wherein $h_i$ is the absolute value of the determinant of the linking matrix of the graph obtained by deleting a terminal vertex on the $i$-th leg from the decomposition $\Gamma$ of $\Sigma(b_1, b_2, b_3)$, and $s$ is the number of vertices in $\Gamma$. 
\end{prop}

\begin{proof}
Let $Y = \Sigma(b_1, b_2, b_3)$. Using the notation from \hyperref[corrected-brie]{Theorem \ref{corrected-brie}} (and that $p:= b_1b_2b_3$ for later bits in the proof) we know that $$\widehat{Z}_0(Y;q) = q^\xi \cdot \left(\widetilde{\Psi}_{b_1b_2b_3}^{(\alpha_1) - (\alpha_2) - (\alpha_3) + (\alpha_4)}(q)\right)$$ where $\xi$ was defined in \hyperref[corrected-brie]{Theorem \ref{corrected-brie}}. Then using the fact that in equation \ref{theta-func-1} we defined that  \begin{equation}
    \widetilde{\Psi}^{(a)}_p(q) := \sum_{n=0}^\infty \psi^{(a)}_{2p}(n) q^{\frac{n^2}{4p}} \in q^{\frac{a^2}{4p}}\mathbb{Z}[[q]]
\end{equation}    we see that $\widetilde{\Psi}_{b_1b_2b_3}^{(\alpha_i)}(q) = q^{\frac{\alpha_i^2}{4p}}P_i(q)$ where $P_i(q) \in \mathbb{Z}[[q]]$ for $1 \leq i \leq 4$. Thus, by making use of how we defined $\widetilde{\Psi}_{b_1b_2b_3}^{(\alpha_1) - (\alpha_2) - (\alpha_3) + (\alpha_4)}(q)$ as a linear combination, we can factorize $\widehat{Z}_0(Y;q)$ into the form:
\begin{align}
    \widehat{Z}_0(Y;q) &= q^\xi \cdot \left(\widetilde{\Psi}_{b_1b_2b_3}^{(\alpha_1)}(q) - \widetilde{\Psi}_{b_1b_2b_3}^{(\alpha_2)}(q) - \widetilde{\Psi}_{b_1b_2b_3}^{(\alpha_3)}(q)+ \widetilde{\Psi}_{b_1b_2b_3}^{(\alpha_4)}(q)\right) \\
    &= q^\xi\left(q^{\frac{\alpha_1^2}{4p}}P_1(q) - q^{\frac{\alpha_2^2}{4p}}P_2(q) - q^{\frac{\alpha_3^2}{4p}}P_3(q) + q^{\frac{\alpha_4^2}{4p}}P_4(q)\right) \\
    &= q^{\xi + \frac{\alpha_1^2}{4p}}\left(P_1(q) - q^{\frac{\alpha_2^2-\alpha_1^2}{4p}}P_2(q) - q^{\frac{\alpha_3^2-\alpha_1^2}{4p}}P_3(q) + q^{\frac{\alpha_4^2-\alpha_1^2}{4p}}P_4(q)\right). \label{fin-form}
\end{align}
Note that \hyperref[brie-pow-int]{Lemma \ref{brie-pow-int}} part (ii) implies that $\frac{\alpha_i^2-\alpha_1^2}{4p} \in \mathbb{Z}$ for $1\leq i \leq 4$. Moreover \hyperref[brie-pow-int]{Lemma \ref{brie-pow-int}} part (i) implies that $\frac{\alpha_i^2-\alpha_1^2}{4p} \in \mathbb{N} \cup \{0\}$ for $1\leq i \leq 4$. Thus in particular we see that $q^{\frac{\alpha_i^2-\alpha_1^2}{4p}}P_i(q) \in \mathbb{Z}[[q]]$ for $1 \leq i \leq 4$. Moreover since $\widetilde{\Psi}_{b_1b_2b_3}^{(\alpha_1)}(q) =  q^{\frac{\alpha_1^2}{4p}}P_1(q)$ we then see by expanding the definition of $\widetilde{\Psi}_{b_1b_2b_3}^{(\alpha_1)}(q)$ that \begin{align*}
    q^{\frac{\alpha_1^2}{4p}}P_1(q) &= \widetilde{\Psi}_{b_1b_2b_3}^{(\alpha_1)}(q)\\
    &= q^{\frac{\alpha_1^2}{4p}} + \sum_{m \geq 1}^\infty q^{\frac{\alpha_1^2}{4p}+\alpha_im+m^2p} - q^{\frac{\alpha_1^2}{4p}-\alpha_im+m^2p} \\
    &= q^{\frac{\alpha_1^2}{4p}}\left(1 + \sum_{m \geq 1}^\infty q^{\alpha_im+m^2p} - q^{-\alpha_im+m^2p}\right)
\end{align*} which implies that $P_1(q) = 1 + \sum_{m \geq 1}^\infty q^{\alpha_im+m^2p} - q^{-\alpha_im+m^2p}$ (to see this simply multiply\footnote{This multiplication is well-defined in $\mathbf{k}$, the Novikov field} both sides of the above equation by $q^{-\frac{\alpha_1^2}{4p}}$) and importantly from this we deduce that the leading term of $P_1(q)$ is $1$. From this we can apply \hyperref[read-delta]{Lemma \ref{read-delta}} to the form of $\widehat{Z}_0(Y;q)$ obtained in equation \ref{fin-form} to see that  $$\Delta_0(Y) = \xi +\frac{\alpha_1^2}{4p}.$$ 
\end{proof}

\subsection{Computing \texorpdfstring{$\widehat{Z}_0$}{Z-hat} and \texorpdfstring{$\Delta_0$}{Delta-0} for Brieskorn Spheres}\label{process-to-compute}
Suppose a Brieskorn sphere $\Sigma(b_1, b_2, b_3)$ is given such that $(b_1, b_2, b_3) \neq (2, 3, 5)$ and one wants to compute $\widehat{Z}_0$ and $\Delta_0$ for it. Then one needs to do the following.

\begin{enumerate}[(i)]
    \item First find integers $b < 0$ and $a_1, a_2, a_3 >0$ which satisfy $b_1b_2b_3 \cdot b + b_2b_3a_1 + b_1b_3a_2 + b_1b_2a_3 = -1$ and  produce the plumbing graph $\Gamma$ of $\Sigma(b_1, b_2, b_3)$.
    \item Then delete the terminal vertices to produce new graphs $\Gamma_i$ with linking matrices $M_i$ for $1 \leq i \leq 3$ and compute $h_i := |\operatorname{det}(M_i)|$. 
    \item Then compute $\alpha_i$ (as defined in Proposition \ref{new-brieskorn-form}) for $1 \leq i \leq 4$.
    \item Finally, use all of this data as input to Proposition \ref{new-brieskorn-form} to compute $\widehat{Z}_0(\Sigma(b_1, b_2, b_3))$ and $\Delta_0(\Sigma(b_1, b_2, b_3))$.
\end{enumerate}

This process can be viewed as an algorithm and thus be coded into a program. Let's now consider an example to see how to use this proposition to compute the $\widehat{Z}_0$ and $\Delta_0$ invariants for Brieskorn spheres in practice.

\begin{exmp}\label{sigma-2-9-11-inv}
Consider $\Sigma(2, 9, 11)$. The plumbing description for this Brieskorn sphere was computed in \hyperref[sigma-2-9-11]{Example \ref{sigma-2-9-11}} and is depicted below:
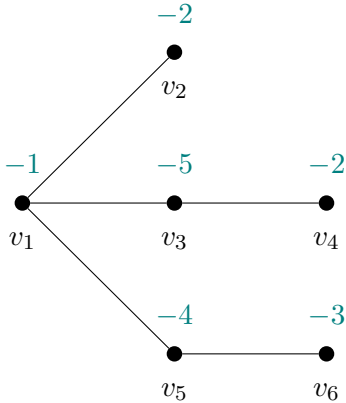
\begin{figure}[H]
\centering

    \begin{tikzpicture}

        \node (0A) at (-4, -0.5) {$v_1$};
        \node (0B) at (-2, 1.5) {$v_2$};
        \node (0C) at (-2, -0.5) {$v_3$};
        \node (0D) at (-2, -2.5) {$v_5$};
        \node (0E) at (0, -0.5) {$v_4$};
        \node (0F) at (0, -2.5) {$v_6$};

        \node (0A) at (-4, 0.5) [text=teal]{$-1$};
        \node (0B) at (-2, 2.5) [text=teal]{$-2$};
        \node (0C) at (-2, 0.5) [text=teal]{$-5$};
        \node (0D) at (-2, -1.5) [text=teal]{$-4$};
        \node (0E) at (0, 0.5) [text=teal]{$-2$};
        \node (0F) at (0, -1.5) [text=teal]{$-3$};
        
        \begin{scope}[every node/.style={circle, thick, draw, fill=black, inner sep =0pt, minimum size=5pt}]
        \node (1) at (-4, 0) {};
        \node (2A) at (-2, 2) {};
         \node (2B) at (-2, 0) {};
          \node (2C) at (-2, -2) {};
          \node (3A) at (0, 0) {};
          \node (3B) at (0, -2) {};
        \end{scope}

        \draw (1) -- (2A);
        \draw (1) -- (2B);
         \draw (1) -- (2C);
         \draw (2B) -- (3A);
           \draw (2C) -- (3B);
    \end{tikzpicture}

\caption{Plumbing graph for $\Sigma(2, 9, 11)$}
\end{figure}
We first notice that $s = |\operatorname{Vert}(\Gamma)| = 6$. The linking matrix is given by: 

\begin{equation*}
        M = \begin{bNiceArray}[first-row]{cccccc} 
        v_1 & v_2 & v_3 & v_4 & v_5 & v_6\\
        -1 & 1 & 1 & 0 & 1 & 0 \\
1 & -2 & 0 & 0 & 0 & 0 \\
1 & 0 & -5 & 1 & 0 & 0 \\
0 & 0 & 1 & -2 & 0 & 0 \\
1 & 0 & 0 & 0 & -4 & 1 \\
0 & 0 & 0 & 0 & 1 & -3
\end{bNiceArray}
\end{equation*}

We find that $\operatorname{Tr}(M) = -17$. If we delete the terminal vertex $v_2$ on the first leg we obtain the linking matrix $M_1$ for the corresponding plumbing graph $\Gamma_1$. Deleting the terminal vertex $v_4$ on the second leg yields the linking matrix $M_2$ for the corresponding plumbing graph $\Gamma_2$. Deleting the terminal vertex $v_6$ on the third leg yields the linking matrix $M_3$ for the corresponding plumbing graph $\Gamma_3$.

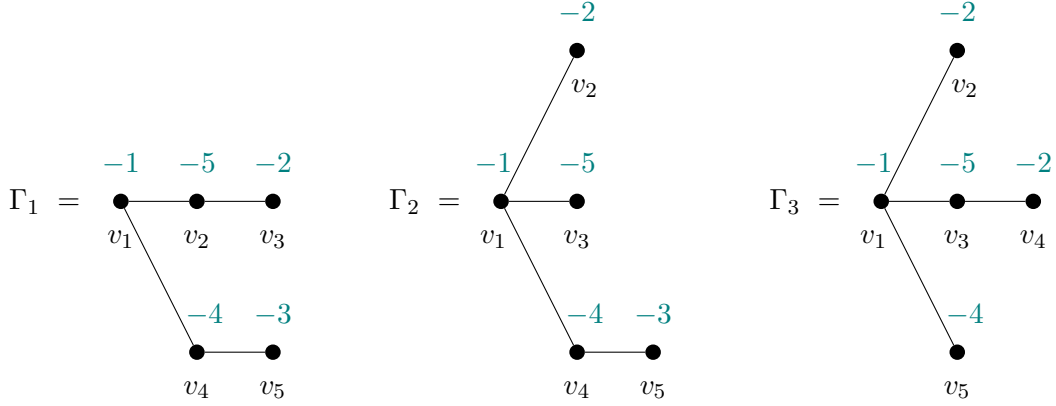
\begin{figure}[H]
    \centering

    \begin{tikzpicture}
        \node (0AA) at (-5, 0) {$\Gamma_1 \ =$};
        \node (0A) at (-4, -0.5) {$v_1$};
        \node (0C) at (-3, -0.5) {$v_2$};
        \node (0D) at (-3, -2.5) {$v_4$};
        \node (0E) at (-2, -0.5) {$v_3$};
        \node (0F) at (-2, -2.5) {$v_5$};

        \node (0A) at (-4, 0.5) [text=teal]{$-1$};
        \node (0C) at (-3, 0.5) [text=teal]{$-5$};
        \node (0D) at (-2.9, -1.5) [text=teal]{$-4$};
        \node (0E) at (-2, 0.5) [text=teal]{$-2$};
        \node (0F) at (-2, -1.5) [text=teal]{$-3$};
        
        \begin{scope}[every node/.style={circle, thick, draw, fill=black, inner sep =0pt, minimum size=5pt}]
        \node (1) at (-4, 0) {};
         \node (2B) at (-3, 0) {};
          \node (2C) at (-3, -2) {};
          \node (3A) at (-2, 0) {};
          \node (3B) at (-2, -2) {};
        \end{scope}

        \draw (1) -- (2B);
         \draw (1) -- (2C);
         \draw (2B) -- (3A);
           \draw (2C) -- (3B);


     \node (0AA) at (0, 0) {$\Gamma_2 \ =$};
        \node (0A) at (0.9, -0.5) {$v_1$};
        \node (0B) at (2.1,  1.5) {$v_2$};
        \node (0C) at (2, -0.5) {$v_3$};
        \node (0D) at (2, -2.5) {$v_4$};
        \node (0F) at (3, -2.5) {$v_5$};

        \node (0A) at (0.9, 0.5) [text=teal]{$-1$};
         \node (0B) at (2, 2.5) [text=teal]{$-2$};
        \node (0C) at (2, 0.5) [text=teal]{$-5$};
        \node (0D) at (2.1, -1.5) [text=teal]{$-4$};
        \node (0F) at (3, -1.5) [text=teal]{$-3$};
        
        \begin{scope}[every node/.style={circle, thick, draw, fill=black, inner sep =0pt, minimum size=5pt}]
        \node (1) at (1, 0) {};
        \node (2A) at (2, 2) {};
         \node (2B) at (2, 0) {};
          \node (2C) at (2, -2) {};
          \node (3B) at (3, -2) {};
        \end{scope}

         \draw (1) -- (2A);
        \draw (1) -- (2B);
         \draw (1) -- (2C);
           \draw (2C) -- (3B);


     \node (0AA) at (5, 0) {$\Gamma_3 \ =$};
        \node (0A) at (5.9, -0.5) {$v_1$};
        \node (0B) at (7.1,  1.5) {$v_2$};
        \node (0C) at (7, -0.5) {$v_3$};
        \node (0D) at (7, -2.5) {$v_5$};
        \node (0E) at (8, -0.5) {$v_4$};

        \node (0A) at (5.9, 0.5) [text=teal]{$-1$};
        \node (0B) at (7, 2.5) [text=teal]{$-2$};
        \node (0C) at (7, 0.5) [text=teal]{$-5$};
        \node (0D) at (7.1, -1.5) [text=teal]{$-4$};
        \node (0E) at (8, 0.5) [text=teal]{$-2$};

        \begin{scope}[every node/.style={circle, thick, draw, fill=black, inner sep =0pt, minimum size=5pt}]
            \node (1) at (6, 0) {};
            \node (2A) at (7, 2) {};
            \node (2B) at (7, 0) {};
            \node (2C) at (7, -2) {};
            \node (3A) at (8, 0) {};
        \end{scope}

        \draw (1) -- (2A);
        \draw (1) -- (2B);
        \draw (1) -- (2C);
        \draw (2B) -- (3A);

    \end{tikzpicture}

    \caption{Plumbing description of the graphs $\Gamma_i$ for $1 \leq i \leq 3$. }
    
\end{figure}

The linking matrices $M_i$ corresponding to $\Gamma_i$ for $1 \leq i \leq 3$ are collected below:
\begin{align*}
        M_1 = \begin{bNiceArray}[first-row]{ccccc} 
        v_1 & v_2 & v_3 & v_4 & v_5 \\
        -1 &  1 & 0 & 1 & 0 \\
1 &  -5 & 1 & 0 & 0 \\
0 &  1 & -2 & 0 & 0 \\
1 &  0 & 0 & -4 & 1 \\
0 &  0 & 0 & 1 & -3
\end{bNiceArray}  M_2 = \begin{bNiceArray}[first-row]{ccccc} 
        v_1 & v_2 & v_3 & v_4 & v_5 \\
         -1 & 1 & 1 &  1 & 0 \\
1 & -2 & 0 &  0 & 0 \\
1 & 0 & -5 &  0 & 0 \\
1 & 0 & 0 &  -4 & 1 \\
0 & 0 & 0 &  1 & -3
\end{bNiceArray}  M_3 = \begin{bNiceArray}[first-row]{ccccc} 
        v_1 & v_2 & v_3 & v_4 & v_5 \\
        -1 & 1 & 1 & 0 & 1 \\
1 & -2 & 0 & 0 & 0  \\
1 & 0 & -5 & 1 & 0  \\
0 & 0 & 1 & -2 & 0  \\
1 & 0 & 0 & 0 & -4 
\end{bNiceArray}
\end{align*}

Now one needs to compute $h_i = |\det(M_i)|$ for $1 \leq i \leq 3$. One finds that $h_1 = 50$, $h_2 = 3$ and $h_3 = 2$. Letting $b_1 = 2, b_2 = 9, b_3 = 11$ and recalling that $p = b_1b_2b_3 = 198$ we  further calculate that $\alpha_1 = 59$, $\alpha_2 = 95$, $\alpha_3= 103$, $\alpha_4 = 139$ and moreover, $\frac{\alpha_2^2 - \alpha_1^2}{4p} = 7$, $\frac{\alpha_3^2 - \alpha_1^2}{4p} = 9$, $\frac{\alpha_4^2 - \alpha_1^2}{4p} = 20$. Then by substituting in the values we just found in the formulas we derived for $\Delta_0$ and $\widehat{Z}_0$ for Brieskorn spheres in Proposition \ref{new-brieskorn-form}, we find that $$\Delta_0(\Sigma(2, 9, 11)) = \frac{9}{2}$$ and
$$\widehat{Z}_0(\Sigma(2, 9, 11);q) =  q^{\frac{9}{2}}\cdot \left[q^{-\frac{59^2}{198}}\left(\widetilde{\Psi}_{198}^{(59) - (95) - (103) + (139)}(q)\right)\right].$$
\end{exmp}

\begin{exmp}\label{sigma-3-7-8-inv}
Consider the Brieskorn sphere $\Sigma(3, 7, 8)$. The integers $b = -1$, $a_1 = 1, a_2 = 2$ and $a_3 = 3$ satisfy equation \ref{brieskorn-eq} with $b_1 = 3, b_2 = 7$ and $b_3=8$. One can then compute the continued fractions $\frac{b_1}{a_1} = 3 = [3]$, $\frac{b_2}{a_2} = \frac{7}{2} = [4,2]$ and $\frac{b_3}{a_3} = \frac{8}{3} = [3,3]$ to produce the following plumbing graph for $\Sigma(3, 7, 8)$. 
\begin{figure}[H]
\centering

 \begin{tikzpicture}

        \node (0A) at (-4, -0.5) {$v_1$};
        \node (0B) at (-2, 1.5) {$v_2$};
        \node (0C) at (-2, -0.5) {$v_3$};
        \node (0D) at (-2, -2.5) {$v_5$};
        \node (0E) at (0, -0.5) {$v_4$};
        \node (0F) at (0, -2.5) {$v_6$};

        \node (0A) at (-4, 0.5) [text=teal]{$-1$};
        \node (0B) at (-2, 2.5) [text=teal]{$-3$};
        \node (0C) at (-2, 0.5) [text=teal]{$-4$};
        \node (0D) at (-2, -1.5) [text=teal]{$-3$};
        \node (0E) at (0, 0.5) [text=teal]{$-2$};
        \node (0F) at (0, -1.5) [text=teal]{$-3$};
        
        \begin{scope}[every node/.style={circle, thick, draw, fill=black, inner sep =0pt, minimum size=5pt}]
        \node (1) at (-4, 0) {};
        \node (2A) at (-2, 2) {};
         \node (2B) at (-2, 0) {};
          \node (2C) at (-2, -2) {};
          \node (3A) at (0, 0) {};
          \node (3B) at (0, -2) {};
        \end{scope}

        \draw (1) -- (2A);
        \draw (1) -- (2B);
         \draw (1) -- (2C);
         \draw (2B) -- (3A);
           \draw (2C) -- (3B);
    \end{tikzpicture}

\caption{Plumbing graph for $\Sigma(3, 7, 8)$}
\end{figure}
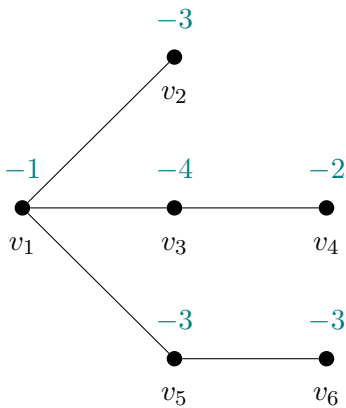
We leave it as an exercise to the reader to use this plumbing description to find $\Delta_0$ and $\widehat{Z}_0$ for $\Sigma(3, 7, 8)$ from Proposition \ref{new-brieskorn-form}. One finds that $$\Delta_0(\Sigma(3, 7, 8)) =  \frac{13}{2}$$ and further that
$$\widehat{Z}_0(\Sigma(3, 7, 8);q) =  q^{\frac{13}{2}}\cdot \left[q^{-\frac{67^2}{168}}\left(\widetilde{\Psi}_{168}^{(67) - (109) - (115) + (157)}(q)\right)\right].$$
\end{exmp}

\section{\texorpdfstring{$\Delta_0$}{Delta} and Homology cobordism}

We begin with a brief review of homology cobordism, the interested reader may consult either \cite{saveliev2011lectures}, \cite{gamkrelidze2013invariants}, \cite{_avk_2023} or \cite{manolesc_cob}. Let $\Sigma_0$ and $\Sigma_1$ be two oriented integral homology $3$-spheres. We call $\Sigma_0$ and $\Sigma_1$, \textit{homology cobordant} if there exists a smooth, compact, oriented $4$-manifold $W$ with boundary $\partial W = -\Sigma_0 \sqcup \Sigma_1$ such that inclusions $\iota_0 : \Sigma_0 \hookrightarrow W$ and $\iota_1 : \Sigma_1 \hookrightarrow W$ induce isomorphisms $(\iota_0)_* : H_i(\Sigma_0) \to H_i(W)$ and $(\iota_1)_* : H_i(\Sigma_1) \to H_i(W)$ for all $i \geq 0$. Homology cobordism is an equivalence relation on the class of all oriented integral homology $3$-spheres.

We define the \textit{($3$-dimensional) homology cobordism group}, $\Theta^3$, to be the class of oriented integral smooth homology $3$-spheres modulo the equivalence relation of homology cobordism. An addition operation is given by $[M] + [N] := [M \# N]$ wherein $[M], [N]$ denote the homology cobordism classes of $M$ and $N$. The identity of the group is $[S^3]$ and for every $[M]$ an additive inverse is given by $[-M]$. Accordingly, we say that an oriented integral smooth homology $3$-sphere is \textit{homology cobordant to zero} if it is homology cobordant to $S^3$.  

A \textit{homology cobordism invariant} is a function $f : \Theta^3 \to X$ where $X$ is any set. In particular this means that if $M$ and $N$ are two oriented homology $3$-spheres which are homology cobordant, then $f(M) = f(N)$. Two examples of homology cobordism invariants are given by the Rokhlin invariant $\mu : \Theta^3 \to \mathbb{Z}/2$ and the correction terms to Heegaard Floer homology $d : \Theta^3 \to \mathbb{Z}$ which we'll see later. A question was raised in \cite{cobordism} which asked whether the $\Delta_0$ invariants were homology cobordism invariants. In this section we will show that this is not the case. 

In searching for a counterexample, we will make use of the following theorem which appears in a different (but equivalent) form in \cite{casson-harer} and \cite[Problem 4.2 - p. 183]{kirby-problems}.

\begin{thm}\label{hom-cob-brie-sphere-thm}

The following families of Brieskorn spheres are all homology cobordant to $S^3$. 
\begin{enumerate}[(i)]
    \item $\Sigma(p, pq-1, pq+1)$ for $p$ even and $q$ odd,
    \item $\Sigma(p, pq+1, pq+2)$ for $p$ odd and any $q$.
\end{enumerate} 
\end{thm}

\begin{proof}
By \cite{casson-harer} and \cite[Problem 4.2 - p. 183]{kirby-problems} manifolds from the families (i) and (ii) bound smooth contractible $4$-manifolds. Let $\Sigma$ denote any manifold from either of the families (i) or (ii) above and let $W$ denote the smooth contractible $4$-manifold $W$ which bounds $\Sigma$. Then since $W$ is a contractible $4$ -manifold with non-empty boundary which is a homology $3$-sphere, $W$ is compact and orientable. It is known (see \cite[Proposition 2.3.4]{fushida-hardy}) that a smooth oriented integral homology $3$-sphere $M$ bounds a smooth compact orientable $4$-manifold $W$ with $H_i(W) = 0$ for all $i \geq 0$, if and only if $M$ is homology cobordant to $S^3$. From this above characterization, we thus conclude that $\Sigma$ is homology cobordant to $S^3$.
\end{proof}

\begin{cor}\label{2-9-11-3-7-8-hom-cob-zero}
The Brieskorn manifolds $\Sigma(2, 9, 11)$ and $\Sigma(3, 7, 8)$ are both homology cobordant to $S^3$. 
\end{cor}

\begin{proof}
From Theorem \ref{hom-cob-brie-sphere-thm}, in (i) take $p=2$ and $q=5$ and in (ii) take $p=3$ and $q=2$. 
\end{proof}

\begin{thm}\label{delta-hom-cob-not}\label{z-hat-cob-not}
$\widehat{Z}_a$ and $\Delta_a$ are neither homology cobordism invariants, nor cobordism invariants.
\end{thm}

\begin{proof}
We saw in Corollary \ref{2-9-11-3-7-8-hom-cob-zero} that $\Sigma(2, 9, 11)$ is homology cobordant to $S^3$. Moreover since any homology cobordism between two manifolds is by definition a cobordism between them, $\Sigma(2, 9, 11)$ is also cobordant to $S^3$.

If in general $\widehat{Z}_a$ and $\Delta_a$ were homology cobordism invariants, then we would have in this particular case that $\widehat{Z}_0(\Sigma(2, 9, 11);q) = \widehat{Z}_0(S^3;q)$ and also that $\Delta_0(\Sigma(2, 9, 11)) = \Delta_0(S^3)$. However in \hyperref[sigma-2-9-11-inv]{Example \ref{sigma-2-9-11-inv}} and \hyperref[delta-s3]{Example \ref{delta-s3}} respectively we computed that
\begin{align*}
    \widehat{Z}_0(\Sigma(2, 9, 11);q) &= q^{\frac{9}{2}}(1-q^7-q^{9}+q^{20}+q^{79}+\cdots) \\
    \widehat{Z}_0(S^3;q) &= q^{-\frac{1}{2}}(2q-2) \\
    \Delta_0(\Sigma(2, 9, 11)) &= \frac{9}{2}\\
    \Delta_0(S^3) &= -\frac{1}{2}.
\end{align*} 
Thus the $\widehat{Z}_a$ and $\Delta_a$ invariants are not invariants of homology cobordism. Similarly since $\Sigma(2, 9, 11)$ is also cobordant to $S^3$ neither the $\widehat{Z}_a$ nor $\Delta_a$ invariants are invariants of cobordism.
\end{proof}

\begin{rem}
In the above we can use other Brieskorn manifolds $\Sigma(b_1,b_2, b_3)$ which appeared in \hyperref[hom-cob-brie-sphere-thm]{Theorem \ref{hom-cob-brie-sphere-thm}}.
\end{rem}

\section{\texorpdfstring{$\Delta_a$}{Delta} and correction terms, \texorpdfstring{$d$}{d}, to Heegaard Floer Homology}

We recap Heegard Floer homology closely following \cite{ozsvath2004introduction} and \cite{manolesc_cob}. 

\subsection{Heegard Floer Homology}

Let $Y$ be a closed orientable $3$-manifold with Heegard decomposition $Y = U_0 \cup_\Sigma U_1$. From the Heegard decomposition we can produce a Heegard diagram $(\Sigma_g, \alpha_1, \dots, \alpha_g, \beta_1, \dots, \beta_g)$ which is sometimes written as $(\Sigma_g, \mathbf{\alpha}, \mathbf{\beta})$. We then define $\operatorname{Sym}(\Sigma_g) = \left(\prod_{i=1}^g\Sigma_g\right)/\mathfrak{S}_g $ where $\mathfrak{S}_g$ is the symmetric group on $g$ letters which acts on $\left(\prod_{i=1}^g\Sigma_g\right)$ in the following way: $\left(\sigma, (x_1, \dots, x_g)\right) \mapsto (x_{\sigma(1)}, \dots, x_{\sigma(g)})$ for $\sigma \in \Sigma_g$ and $(x_1, \dots, x_g) \in \prod_{i=1}^g\Sigma_g$  (see \cite[Section 4]{ozsvath2004introduction}). As similarly stated in \cite[Section 4.1]{ozsvath2004introduction}, the "attaching circles" $\alpha_i$ and $\beta_i$ produce tori $$\mathbb{T}_\alpha = \alpha_1 \times \cdots \times \alpha_g \text{ and } \mathbb{T}_\beta = \beta_1 \times \cdots \times \beta_g$$ which are subsets of $\operatorname{Sym}(\Sigma_g)$. A certain map $s_z : \mathbb{T}_\alpha \cap \mathbb{T}_\beta \to \operatorname{Spin}^c(Y)$ is then defined (see \cite[Section 6.2]{ozsvath2004introduction}) from which we can define groups:
\begin{enumerate}
    \item $CF^\infty(\alpha, \beta, a)$ to be the free abelian group generated by the pairs $[x, i]$ where the $x \in \mathbb{T}_\alpha \cap \mathbb{T}_\beta$ with $s_z(x) = a$ and $i \in \mathbb{Z}$. 
    \item $CF^{-}(\alpha, \beta, a)$ to be the free abelian group generated by the pairs $[x, i]$ where the $x \in \mathbb{T}_\alpha \cap \mathbb{T}_\beta$ with $s_z(x) = a$ and $i <0$ where $i$ is an integer.
    \item $CF^{+}(\alpha, \beta, a) := CF^\infty(\alpha, \beta, a)/CF^{-}(\alpha, \beta, a)$.
\end{enumerate}

One can assign a (relative) grading to each of these groups which then allow us to define chain complexes from them. From the resulting chain complexes we get homology groups $HF^{-}(Y, a), HF^+(Y, a)$ and  $HF^{\infty}(Y, a)$ as is usually done in algebraic topology. 

According to \cite[Theorem 7.1]{os-hol-triangle}, if $Y$ is additionally a rational homology sphere, we have that $HF^{-}(Y, a), HF^+(Y, a), HF^{\infty}(Y, a)$ are $\mathbb{Q}$-graded. This means, for example in the case of $HF^+(Y, a)$ , that $HF^+(Y, a) = \bigoplus_{\omega \in \mathbb{Q}} HF_\omega^+(Y, a).$ We have a similar direct sum decomposition for $HF^{\infty}(Y, a)$. 

Furthermore for each $\omega \in \mathbb{Q}$ there is a family of homomorphisms $HF_\omega^+(Y, a) \to HF^{\infty}_\omega(Y, a)$ which come from a long exact sequence (see \cite[Section 2, eq. 3]{oz_szabo_main}) $$\cdots HF^{-}(Y, a) \to HF^+(Y, a) \to HF^{\infty}(Y, a) \to \cdots.$$ 

The following is a slightly expanded version of \cite[Definition 4.1]{oz_szabo_main}.
\begin{defn}
Let $Y$ be an oriented rational homology $3$-sphere and $a \in \operatorname{Spin}^c(Y)$.  We define a rational number $d(Y, a) \in \mathbb{Q}$, called the \textit{correction term}, to be the minimal $\omega \in \mathbb{Q}$ such that an element $x$ in the image of the homomorphism $HF_\omega^+(Y, a) \to HF^{\infty}_\omega(Y, a)$ is non-torsion, i.e. $x^n \neq 0$ for any integer $n$.
\end{defn}

\begin{rem}
In the case when $Y$ is an oriented $3$-manifold with $b_1(Y) = 0$, then there is a unique $\operatorname{Spin}^c$ structure on $Y$ and we simply write $d(Y)$ instead of $d(Y, a)$. 
\end{rem}

\begin{thm}[Properties of correction terms]\label{correction-term-props}
\leavevmode
\begin{enumerate}[(i)]
    \item The correction terms give group homomorphisms:
    \begin{enumerate}[(a)]
        \item $d : \theta^c \to \mathbb{Q}$ where $\theta^c$ is the $\operatorname{Spin}^c$ homology cobordism group (\cite[Theorem 1.2]{oz_szabo_main})
        \item $d : \Theta^3 \to \mathbb{Z}$ where $\Theta^3$ is the homology cobordism group (\cite[p. 9]{manolesc_cob})
    \end{enumerate}
    \item Let $(Y, a)$ be a rational homology $3$-sphere with $\operatorname{Spin}^c$ structure $a$. Then it follows that $d(-Y, a) = - d(Y, a)$ where $-Y$ is the orientation reversal of $Y$. (\cite[Proposition 4.2]{oz_szabo_main})
    \item Let $(Y, a)$ be a rational homology $3$-sphere with $\operatorname{Spin}^c$ structure $a$. Then it follows that $d(Y, a) = d(Y, \overline{a})$ where $\overline{a}$ is the conjugation of $a$. (\cite[Proposition 4.2]{oz_szabo_main})
\end{enumerate}
\end{thm}

The following is a corollary to part (i) of Theorem \ref{correction-term-props} above.

\begin{cor}
The correction term for $S^3$ is zero, i.e. $d(S^3) = 0$
\end{cor}

\begin{proof}
The homology cobordism class of $S^3$ is the identity of the group $\Theta^3$. Thus it must be mapped to $0 \in \mathbb{Z}$ since $d$ is a group homomorphism. 
\end{proof}

\subsection{Connections between \texorpdfstring{$\Delta_a$}{Delta} and \texorpdfstring{$d$}{d}}

Notice that for an oriented rational homology $3$-sphere $Y$ with $\operatorname{Spin}^c$ structure $a$, we have the following connections, outlined in \cite{cobordism}, between  $\Delta_a(Y)$ and $d(Y, a)$:
\begin{enumerate}[(i)]
    \item $\Delta_a(Y)$ and $d(Y, a)$ are both labelled by $\operatorname{Spin}^c$ structures on $Y$.
    \item $\Delta_a(Y)$ and $d(Y, a)$ both remain unchanged under conjugation of $\operatorname{Spin}^c$ structures.
\end{enumerate}
Therefore it is natural to ask the question: \textit{"Just how closely related are $\Delta_a(Y)$ and $d(Y, a)$?"}. In particular $d(Y, a)$ is a homology cobordism invariant, so a further question one could ask is: \textit{"Is $\Delta_a(Y)$ also a homology cobordism invariant"?} To partially answer these questions, we have the following proposition from \cite{cobordism}, which we state below in the case of \textit{almost rational graphs} (see \cite{nemethi_ozvath_szabo}).

\begin{prop}\label{d-conj}
If $Y = Y(\Gamma)$ is a negative definite plumbed manifold arising from an almost rational graph $\Gamma$ then $\Delta_a(Y) = \frac{1}{2} - d(Y, a) \mod 1.$
\end{prop}

\begin{rem}
The result that $\Delta_a(Y) = \frac{1}{2} - d(Y, a) \mod 1$ could equivalently be rewritten as as $\Delta_a(Y) = n+ \frac{1}{2} - d(Y, a)$ for some $n \in \mathbb{Z}$.
\end{rem}

The following corollary follows immediately from the above proposition and the fact that the correction terms, $d$, are homology cobordism invariants.

\begin{cor}
Let $\mathfrak{X}$ be the class of negative definite plumbed manifolds arising from almost rational graphs
The function $f: \mathfrak{X} \to \mathbb{Q}$ defined by $f(Y) = \Delta_a(Y) \mod 1$ is a homology cobordism invariant
\end{cor}

Let us now take a look at an example\footnote{Part of this example, and the corollary above, was inspired by discussion with Mrunmay Jagadale.} to see this phenomenon.

\begin{exmp}
We know from \hyperref[2-9-11-3-7-8-hom-cob-zero]{Corollary \ref{2-9-11-3-7-8-hom-cob-zero}} that $\Sigma(2, 9, 11)$ and $\Sigma(3, 7, 8)$ are homology cobordant to $S^3$. We also saw in Examples \ref{delta-s3}, \ref{sigma-2-9-11-inv} and \ref{sigma-3-7-8-inv} that $\Delta_0(S^3) = \frac{1}{2}$, $\Delta_0(\Sigma(2, 9, 11)) = \frac{9}{2} =  \frac{1}{2} + 4$ and $\Delta_0(\Sigma(3, 7, 8)) = \frac{13}{2} = \frac{1}{2} + 6$. Then notice that 
$\Delta_0(S^3) = \Delta_0(\Sigma(2, 9, 11)) = \Delta_0(\Sigma(3, 7, 8)) \mod 1$ as expected. 
\end{exmp}

\subsection{\texorpdfstring{$\Delta_0$}{Delta} and integral homology spheres}

We now investigate the form of the $\Delta_a$ invariants for negative-definite plumbed manifolds (arising from almost rational graphs) which are also integer homology spheres. Using theory about the correction terms, we find that the $\Delta_0$ invariants for this class of manifolds have a simple form as stated by the following proposition.

\begin{prop}\label{delta-hom-sphere}
    Let $\Gamma$ be an almost rational graph and let $Y := Y(\Gamma)$ be a negative definite plumbed $3$-manifold which is also an integral homology sphere, then $\Delta_0(Y) = \frac{1}{2} \mod 1$
\end{prop}

\begin{proof}
By \hyperref[d-conj]{Proposition \ref{d-conj}} $\Delta_0(Y) = \frac{1}{2} + d(Y) + n$ for some $n \in \mathbb{Z}$. Since $d : \Theta^3 \to \mathbb{Z}$ is a homomorphism we have that $d(Y) \in \mathbb{Z}$ which completes the proof.
\end{proof}

\begin{rem}
    In particular the above proposition gives an explanation for why when $Y$ is a Brieskorn sphere $\Sigma(b_1, b_2, b_3)$ we find that $\Delta_0(Y) = \frac{1}{2} \mod 1$ as observed in all such examples within this text. See also Tables \ref{tab:z-hat-delta-hard} and \ref{tab:further-comp-hom-s3}.
\end{rem}

\subsection{Sharpness of the relation between \texorpdfstring{$\Delta_a$}{Delta invariants} and \texorpdfstring{$d$}{correction terms}}

Suppose now that \hyperref[d-conj]{Proposition \ref{d-conj}} instead stated that $\Delta_a(Y) = \frac{1}{2} - d(Y, a) \mod x$ where $x$ is some integer. Intuitively the possibility of showing that the relation $\Delta_a(Y) = \frac{1}{2} - d(Y, a) \mod x$ holds, for a higher value of $x$ is desirable since this leads to a stronger relation between $\Delta_a(Y)$ and $d(Y, a)$.  It was stated in \cite[p. 25, 26]{cobordism} that the best hope would be to find such a relation $\Delta_a(Y) = \frac{1}{2} - d(Y, a) \mod x$ for $x=2$. This conclusion was based off the examples computed within \cite{cobordism}. Based on the examples in this text, we are able to give an independent proof of this.

\begin{lem}\label{rel-sharp-less2}
Suppose that the relation $\Delta_a(Y) = \frac{1}{2} - d(Y, a) \mod x$ holds for all negative definite plumbed manifolds $Y = Y(\Gamma)$. Then $x \leq 2$.
\end{lem}

\begin{proof}
The relation that  $\Delta_a(Y) = \frac{1}{2} - d(Y, a) \mod x$ can be equivalently restated as \begin{equation}\label{delta-d-relation-pre}
    \Delta_a(Y) = nx + \frac{1}{2} - d(Y, a)
\end{equation} for some $n \in \mathbb{Z}$. Now we have that $\Delta_a(\Sigma(2, 9, 11)) = \frac{9}{2} = \frac{1}{2} + 4$ and $\Delta_a(\Sigma(3, 7, 8)) = \frac{13}{2} = \frac{1}{2} + 6$ since we have that $d(\Sigma(2, 9, 11)) = d(\Sigma(3, 7, 8)) = d(S^3) = 0$ since $\Sigma(2, 9, 11)$ and $\Sigma(3, 7, 8)$ are both homology cobordant to $S^3$. Now by equation \eqref{delta-d-relation-pre} we expect that $4 = nx$ and $6 = mx$ for some $n, m \in \mathbb{Z}$. Thus we see that $x$ is a divisor of both $4$ and $6$. The only possible common divisors are $1$ and $2$ hence $x \in \{1, 2\}$ thus completing the proof.
\end{proof}

\subsection{\texorpdfstring{$\Delta_a$}{Delta-invariants} and \texorpdfstring{$\operatorname{Spin}^c$}{Spin-c} homology cobordism}

Similar to the homology cobordism group $\Theta^3$, there is a $\operatorname{Spin}^c$ homology cobordism group $\theta^3$, (see \cite{ozsvath2002absolutelygradedfloerhomologies}), for which the correction terms provide an invariant of the form $d : \theta^3 \to \mathbb{Q}$. We make the following conjecture.

\begin{conjecture}\label{spinc-hom-cob-not}
$\Delta_a$ is not a $\operatorname{Spin}^c$ homology cobordism invariant.
\end{conjecture}

It's possible that the examples given within our text might be able to be used to verify this conjecture.

\appendix

\section{Comparing \texorpdfstring{$\Delta_0$}{Delta} and \texorpdfstring{$d$}{the correction terms} for some further examples}\label{correction-terms-comparison}

Let us look at some further examples with which we can compare $\Delta_0$ and the correction terms $d$. For many classes of $3$-manifolds there exist techniques which aid the computation of both Heegard Floer homology and the correction terms. Here is one such example that we shall use in the next section. The paper \cite{nem11} gives us the following computational tool.

\begin{prop}\label{d-brie-family-p}
For $\frac{1}{1}$-Dehn surgery on the torus knot $T_{p, p+1}$ embedded in $S^3$, that is for $S^3_{+1}(T_{p, p+1}) = -\Sigma(p, p+1, p(p+1)-1)$ we have $d\left(S^3_{+1}(T_{p, p+1})\right) =- \big \lfloor \frac{p}{2} \big\rfloor\left(\big \lfloor \frac{p}{2} \big\rfloor +1\right).$ 
\end{prop}

Using \hyperref[d-brie-family-p]{Proposition \ref{d-brie-family-p}} we will compute the correction term, $d$, for various manifolds of the form  $-\Sigma(p, p+1, p(p+1)-1)$ and compare the values of $d$ obtained to $\Delta_0$ for these manifolds. These examples are new in that they haven't been explicitly produced in this form before, but are not new in the sense that they say anything new about Proposition \ref{d-conj}. \\ 

There is a subtlety that the reader should be aware of however, as stated in \cite{gukov2020twovariable}, $S^3_{+1}(T_{p, p+1}) = -\Sigma(p, p+1, p(p+1)-1)$ is not a negative definite plumbed manifold, however $\Sigma(p, p+1, p(p+1)-1)$ is a negative definite plumbed manifold, thus we will compare $\Delta_0$ and $d$ for manifolds of the form $\Sigma(p, p+1, p(p+1)-1)$.  To compute $d(\Sigma(p, p+1, p(p+1)-1))$ we simply use the fact that $d(\Sigma(p, p+1, p(p+1)-1)) = -d(-\Sigma(p, p+1, p(p+1)-1))$.\\

We produce the following examples of Brieskorn spheres $\Sigma(p, p+1, p(p+1)-1)$ for which we can compare $\Delta_0$ to $d$. 

\begin{table}[H]
    \centering
    \def\arraystretch{1.5}%
\begin{tabular}{ |c|c|c|c| } 
 \hline
$\Sigma(b_1, b_2, b_3)$ & $\widehat{Z}_0(\Sigma(b_1, b_2, b_3) ;q)$ & $\Delta_0$ & $d$ \\
\hline
$\Sigma(3, 4, 11)$ & $q^{1/2}\left(1-q^5-q^{19}-q^{29}+ \cdots\right)$ & $\frac{1}{2}$ & $2$\\ 
$\Sigma(4,5,19)$ & $q^{37/2}\cdot\left(1-q^{11}-q^{53}-q^{71}+q^{72}+q^{92}+ \cdots\right)$ & $\frac{37}{2}$ & $6$\\ 
$\Sigma(5,6,29)$ & $q^{141/2}\cdot\left(1-q^{19}-q^{111}-q^{139}+q^{140}+ \cdots\right)$ & $\frac{141}{2}$ & $6$\\ 
$\Sigma(6,7,41)$ & $q^{361/2}\cdot\left(1-q^{29}-q^{199}-q^{239}+q^{240} + \cdots\right)$ & $\frac{361}{2}$ & $12$\\
 \hline
\end{tabular}
    \caption{\textsc{Comparisons of $\Delta_0$ and $d$ for various $\Sigma(p, p+1, p(p+1)-1)$}}
    \label{tab:comparison-d}
\end{table}

\section{Some computations of \texorpdfstring{$\widehat{Z}_0$}{Z-hat} and \texorpdfstring{$\Delta_0$}{Delta} for \texorpdfstring{$\Sigma(b_1, b_2, b_3)$}{Brieskorn spheres}}\label{program}

In this final section we record some explicit computations of both $\widehat{Z}_0$ and $\Delta_0$ for various Brieskorn spheres which have not yet appeared in the literature. The tables below highlight the growth/behaviour of the $\Delta_0$ invariants. In particular Table \ref{tab:z-hat-delta-hard} shows computation of $\Delta_0$ for plumbings which a large number of vertices. Furthermore, Table \ref{further-comp-hom-cob} shows computations of $\Delta_0$ for Brieskorn spheres of the form $\Sigma(p, pq-1, pq+1)$ which we know to be homology cobordant to $S^3$.

\begin{table}[H]
    \centering
    \def\arraystretch{1.5}%
\begin{tabular}{ |c|c|c|c| } 
 \hline
$\Sigma(b_1, b_2, b_3)$ & $\widehat{Z}_0(\Sigma(b_1, b_2, b_3) ;q)$ & $\Delta_0$ & No. vertices \\
\hline
$\Sigma(8,35,93)$ & $q^{9045/2}\cdot\left(1-q^{237}-q^{643}+q^{896}-q^{3127}+q^{3434}+ \cdots\right)$ & $\frac{{9045}}{2}$ & $104$\\ 
$\Sigma(17,41,87)$ & $q^{24801/2}\cdot\left(1-q^{639}-q^{1375}+q^{2048}-q^{3439}+q^{4160}+\cdots\right)$ & $\frac{{24801}}{2}$ & $103$\\ 
$\Sigma(17,53,100)$ & $q^{37441/2}\cdot\left(1-q^{831}-q^{1583}+q^{2448}-q^{5147} + \cdots\right)$ & $\frac{{37441}}{2}$&$117$\\ 
$\Sigma(29, 50, 69)$ & $q^{43317/2}\cdot\left(1-q^{1371}-q^{1903}-q^{3331}+q^{3332}+\cdots\right)$ & $\frac{{43317}}{2}$&$119$\\ 
$\Sigma(29, 53, 96)$ & $q^{64617/2}\cdot\left(1-q^{1455}-q^{2659}+q^{4172}-q^{4939}+\cdots\right)$ & $\frac{{64617}}{2}$&$109$\\ 
$\Sigma(31, 61, 63)$ & $q^{52081/2}\cdot\left(1-q^{1799}-q^{1859}-q^{3719}+q^{3720}+\cdots\right)$ & $\frac{{52081}}{2}$&$124$\\ 
$\Sigma(35, 61, 97)$ & $q^{92365/2}\cdot\left(1-q^{2039}-q^{3263}+q^{5372}-q^{5759} + \cdots \right)$ & $\frac{{92365}}{2}$&$117$\\ 
$\Sigma(39, 41, 94)$ & $q^{66265/2}\cdot\left(1-q^{1519}-q^{3533}-q^{3719}+q^{5130}+q^{5320} + \cdots \right)$ & $\frac{{66265}}{2}$&$102$\\ 
$\Sigma(41, 51, 95)$ & $q^{88737/2}\cdot\left(1-q^{1999}-q^{3759}-q^{4699}+q^{5840} + \cdots \right)$ & $\frac{{88737}}{2}$&$109$\\ 
$\Sigma(42, 43, 95)$ & $q^{76141/2}\cdot\left(1-q^{1721}-q^{3853}-q^{3947}+q^{5658}+q^{5754} + \cdots \right)$ & $\frac{{76141}}{2}$&$105$\\ 
 \hline
\end{tabular}
    \caption{\textsc{Finding $\widehat{Z}_0$ and $\Delta_0$ for various $\Sigma(b_1, b_2, b_3)$ with large plumbing descriptions}}
    \label{tab:z-hat-delta-hard}
\end{table}

\subsection{Further computations for \texorpdfstring{$\Sigma(b_1, b_2, b_3)$}{Brieskorn Spheres} homology cobordant to \texorpdfstring{$S^3$}{the three-sphere}.}\label{further-comp-hom-cob}

\begin{table}[H]
    \centering
    \def\arraystretch{1.5}
\begin{tabular}{ |c|c|c| } 
 \hline
$\Sigma(p,pq-1,pq+1)$ & $\widehat{Z}_0(\Sigma(p,pq-1,pq+1);q)$ & $\Delta_0$ \\
\hline
$\Sigma(2,13,15)$ & $q^{25/2}\cdot\left(1-q^{11}-q^{13}+q^{28}-q^{167}+q^{204}+ \cdots\right)$ & $\frac{{25}}{2}$ \\ 
$\Sigma(2, 21, 23)$ & $q^{81/2}\left(1-q^{19}-q^{21}+q^{44} + \cdots \right)$ & $\frac{81}{2}$ \\
$\Sigma(2,81,83)$ & $q^{1521/2}\cdot\left(1-q^{79}-q^{81}+q^{164}-q^{6559}+q^{6800}+\cdots\right)$ & $\frac{{1521}}{2}$ \\ 
$\Sigma(4,11,13)$ & $q^{97/2}\cdot\left(1-q^{29}-q^{35}+q^{72}-q^{119}+q^{170}+q^{180}+ \cdots\right)$ & $\frac{{97}}{2}$ \\ 
$\Sigma(4,59,61)$ & $q^{3697/2}\cdot\left(1-q^{173}-q^{179}+q^{360}-q^{3479}+q^{3770} + \cdots\right)$ & $\frac{{3697}}{2}$ \\ 
$\Sigma(6,17,19)$ & $q^{505/2}\cdot\left(1-q^{79}-q^{89}+q^{180}-q^{287}+q^{400}+\cdots\right)$ & $\frac{{505}}{2}$ \\ 
$\Sigma(6,41,43)$ & $q^{3265/2}\cdot\left(1-q^{199}-q^{209}+q^{420}+q^{1960}-q^{1679} + \cdots\right)$ & $\frac{{3265}}{2}$ \\ 
$\Sigma(8,23,25)$ & $q^{1441/2}\cdot\left(1-q^{153}-q^{167}+q^{336}-q^{527}+q^{726}+\cdots\right)$ & $\frac{{1441}}{2}$\\ 
$\Sigma(8,87,89)$ & $q^{22497/2}\cdot\left(1-q^{601}-q^{615}-q^{7567}+q^{8342}+q^{8360}+\cdots\right)$ & $\frac{{22497}}{2}$\\ 
 \hline
\end{tabular}
    \caption{\label{tab:further-comp-hom-s3} \textsc{Further computations of $\widehat{Z}_0$ and $\Delta_0$ for various $\Sigma(p, pq-1, pq+1)$}}
\end{table}

\printbibliography

\end{document}